\documentclass{amsart}
\usepackage{amssymb}

\usepackage[dvips]{graphicx}
\usepackage[centertags]{amsmath}
\usepackage{amsthm}

\usepackage{epsf}
\usepackage{amsfonts}
\usepackage[mathcal]{eucal}
\usepackage{mathrsfs}
\usepackage{bm}
\usepackage{graphics, labelfig}
\newtheorem{theorem}{Theorem}

\newtheorem{definition}[theorem]{Definition}

\newtheorem{lemma}[theorem]{Lemma}
\newtheorem{sublemma}[theorem]{Sublemma}

\newtheorem{problem}[theorem]{Problem}
\newtheorem{proposition}[theorem]{Proposition}
\newtheorem{remark}[theorem]{Remark}

\def\geq{\geqslant}
\def\leq{\leqslant}

\begin{document}
\title{QUANTUM TRACES IN QUANTUM TEICHM\"{U}LLER THEORY}
\author{Chris Hiatt}
\thanks{This work was partially supported by NSF grant DMS-0103511 at The University of Southern California}

\address {C. Hiatt, Department
of Mathematics and Computer Science,  The University of Texas of
the Permian Basin, Midland, TX, U.S.A.} \email{hiatt\_c@utpb.edu}

\begin{abstract} We prove that for the torus with one hole
and $p \geq 1$ punctures and the sphere with four holes there is a family of quantum trace functions
in the quantum Teichm\"uller space, analog to the non-quantum trace functions in Teichm\"uller space,
 satisfying the properties
proposed by Chekhov and Fock in \cite{CF1}.
\end{abstract}

\date{\today}
\maketitle

\section{Introduction}

The physicists L. Chekhov and V. Fock and, independently, R.
Kashaev introduced a quantization of the Teichm\"{u}ller space as
an approach to quantum gravity in $2+1$ dimensions. A widespread
philosophy in mathematics is that studying a space is the same as
studying the algebra of functions on that space. The quantum
Teichm\"{u}ller space of Chekhov-Fock and Kashaev
$\mathcal{T}_S^q$ is a certain non-commutative deformation of the
algebra of rational functions on the usual Teichm\"{u}ller space
$\mathcal{T}(S)$. Namely, $\mathcal{T}_S^q$ depends on a parameter
$q=e^{i\hbar}$ and converges to the algebra of functions on
$\mathcal{T}(S)$ as $q \rightarrow 1$ or, equivalently as the
Planck constant $\hbar \rightarrow 0$.

    At this point in time, the quantum Teichm\"{u}ller space is only
defined for surfaces with punctures. Namely, the surface $S$ must
be obtained by removing finitely many points from a compact
surface $\bar{S}$ ; and this in such a way that at least one point
is removed from each boundary component and that, when $\partial
\bar{S} = \phi$, at least one point is removed.

There actually are two versions of the quantum Teichm\"{u}ller
space. the ``logarithmic" version is the original version
developed by the physicists \cite{CF1}. The ``exponential" version
was developed by F. Bonahon and X. Liu \cite{Liu1} \cite{BonLiu1}
and is better adapted to mathematics. For instance, the
exponential version has an interesting finite dimensional
representation theory, which turns out to be connected to
hyperbolic geometry \cite{BonLiu1}.

    A simple closed curve $\alpha$ on the surface $S$ determines a
trace function $T_{\alpha} : \mathcal{T}(S) \rightarrow
\mathbb{R}$ defined as follows: If a point of $\mathcal{T}(S)$ is
represented by a group homomorphism $r: \pi_1 (S) \rightarrow
SL_2(\mathbb{R})$ , then $\mathcal{T}_{\alpha} (r)$ is the trace
of $r(\alpha) \in SL_2 (\mathbb{R})$ . Much of the structure of
the Teichm\"{u}ller space $\mathcal{T}(S)$ can be reconstructed
from these trace functions. See \cite{CS1, Gol1, Luo1}.

    In \cite{CF1} Chekhov and Fock proposed the following problem:

\begin{problem}
\label{prob:MainProb}

 For every simple closed curve $\alpha$ on
$S$, there is a quantum analogue  $T_{\alpha}^q$ of the trace
function $T_{\alpha}$ such that:
\begin{enumerate}
    \item $T_{\alpha}^q \in \mathcal{T}_S^q$ is well defined, independent of choice of
    coordinates.
    \item as $q \rightarrow 1$, $T_{\alpha}^q$ converges to the non-quantum
    trace function $T_{\alpha}$ in $\mathcal{T}(S)$
    \item If $\alpha$ and $\beta$ are disjoint, $T_{\alpha}^q$ and
    $T_{\beta}^q$ commute.
    \item If $\alpha$ and $\beta$ meet in one point, and if
    $\alpha \beta$ and $\beta \alpha$ are obtained by resolving
    the intersection point as in Figure 1, then $T_{\alpha}^q
    T_{\beta}^q = q^{\frac{1}{2}} T_{\alpha \beta}^q + q^{-\frac{1}{2}} T_{\beta \alpha}^q$.
    \item If $\alpha$ and $\beta$ meet in two points of opposite algebraic intersection sign, and if
    $\alpha \beta$, $\beta \alpha$, $\gamma_1$, $\gamma_2$, $\gamma_3$, and $\gamma_4$ are obtained by resolving
    the intersection points as in Figure 2, then $T_{\alpha}^q
    T_{\beta}^q = q T_{\alpha \beta}^q + q^{-1} T_{\beta \alpha}^q + T_{\gamma_1}^q T_{\gamma_2}^q + T_{\gamma_3}^q T_{\gamma_4}^q$.
\end{enumerate}

\end{problem}

It can be shown that, if the quantum trace functions
$T_{\alpha}^q$ exist, they are unique by conditions (4) and (5).
Compare for instance \cite{Luo1}.


\begin{figure}[htb]
\centering
 \SetLabels
 ( 0 * .2) $\alpha$\\
  ( .2 * .2) $\beta$\\
  ( .5 * .2) $\alpha\beta$\\
  ( .9 * .06) $\beta\alpha$\\
\endSetLabels
\centerline{\AffixLabels{
\includegraphics{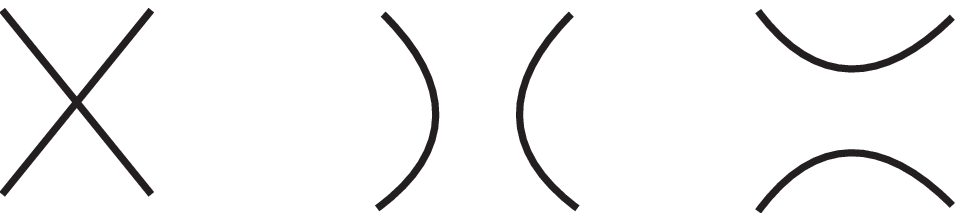}}}
\caption{Resolving Single Crossing} \label{Fig:TorusResolved}
\end{figure}


\begin{figure}[htb]
\centering
 \SetLabels
 ( .11 * .85) $\alpha$\\
  ( .11 * .1) $\beta$\\
  ( .38 * .2) $\alpha\beta$\\
  ( .66 * .2) $\beta\alpha$\\
  ( .9 * .15) $\gamma_1$\\
  ( .9 * .47) $\gamma_2$\\
  ( .9 * .8) $\gamma_3$\\
  ( 1.02 * .5) $\gamma_4$\\
\endSetLabels
\centerline{\AffixLabels{
\includegraphics{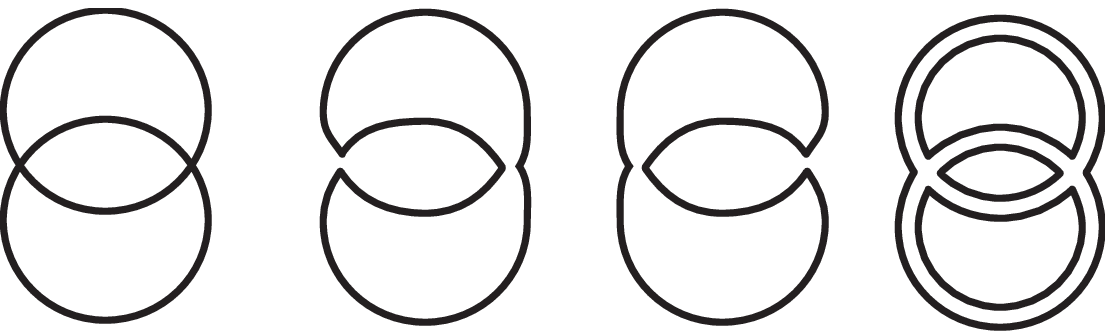}}}
\caption{Resolving Double Crossing}
\end{figure}

\vspace{3mm}

In \cite{CF1}, Chekhov and Fock have verified Problem
~\ref{prob:MainProb} for the once-punctured torus, obtained by
removing one point from the torus, in the case of the logarithmic
model of the quantum Teichm\"{u}ller space.

The exponential model offers some technical challenges, because
certain issues involving square roots have to be resolved to make
sense of Problem~\ref{prob:MainProb}, in particular with respect
to coordinate changes.

The coordinate change isomorphisms introduced by Chekhov-Fock
\cite{CF1}, Kashaev \cite{Kash1}, and Liu \cite{Liu1} satisfy the
following:

\begin{theorem}
\label{thm:CFCoordChanges} " (Chekhov-Fock, Kashaev, Liu) There
exists a family of algebra isomorphisms

\begin{equation}
\Phi_{\lambda \lambda'}:\mathcal{T}_{\lambda'}^{q} \rightarrow
\mathcal{T}_{\lambda}^{q} \nonumber
\end{equation}
indexed by pairs of ideal triangulations $\lambda$, $\lambda'$
of a punctured surface $S$, which satisfy the following conditions:
\begin{enumerate}
\item $\Phi_{\lambda \lambda''} = \Phi_{\lambda \lambda'} \circ
\Phi_{\lambda' \lambda''}$ for any ideal triangulations $\lambda$, $\lambda'$, and
$\lambda''$ of $S$.

\item If $\lambda' = \sigma \lambda$ is obtained by reindexing
$\lambda$ by a permutation $\sigma \in S_n$, then $\Phi_{\lambda
\lambda'}(X_{i}^{'})=X_{\sigma(i)}^{'}$ for any $1 \leq i \leq n$.

\end{enumerate}
\end{theorem}

The first part of this paper is devoted to resolving these
technical issues in the exponential model for the quantum
Teichm\"{u}ller space. This part culminates in the following
theorem.

\begin{theorem}
There is a family of linear maps in the exponential model for the quantum
Teichm\"{u}ller space which
satisfy the conditions of Theorem ~\ref{thm:CFCoordChanges}.

\end{theorem}

The second part of this paper solves Problem~\ref{prob:MainProb} for surfaces
which are at one level of complexity higher that the
once-punctured torus.

We consider the case of the \emph{torus with a wide hole and $p
\geq 1$ punctures}, namely a surface obtained from the compact
surface of genus one with one boundary component by removing $p
\geq 1$ punctures from its boundary, but none from its interior.

\begin{theorem}
\label{thm:MainThmBefore}

If the surface $S$ is a torus with a wide hole and $p \geq 1$
punctures, then there exists a (unique) family of quantum trace
functions as in Problem ~\ref{prob:MainProb}.

\end{theorem}

We then investigate the case of the \emph{sphere with four holes},
where the holes can be either wide or just punctures. Namely, such
a surface is obtained from the compact surface of genus zero with
$k$ boundary components by removing $p$ points from its interior
and at least one point from each boundary component, with $k+p=4$.

\begin{theorem}
\label{thm:tracesinteichBefore}

If $S$ is a sphere with four holes, then there exists a (unique)
family of quantum trace functions as in Problem
~\ref{prob:MainProb}.

\end{theorem}

The overall organization of this paper is as follows:  We
introduce the classical Teichm\"uller Space and the traces in the
non-quantum context. Then we introduce the ``exponential" model of the quantum Teichm\"uller Space and the
analogous quantum traces. We then resolve the technical issues arising from the square roots.
 Finally we prove Theorem
~\ref{thm:MainThmBefore} and Theorem
~\ref{thm:tracesinteichBefore}.

Acknowledgements: I would like to thank Francis Bonahon for his constant help and encouragement. I couldn't have done this paper without his help. I would also like to thank God for giving me the ability to complete this paper.

\section{Ideal Triangulations}

Throughout this paper we will consider an oriented surface S of
finite topological type, where $S = \bar{S} - \{v_1,...,v_p$\} is
obtained by removing $p$ points \{$v_1,...,v_p$\} from a compact
oriented surface \={S} of genus g, with $d \geq 0$ boundary
components. The requirements for S are that $p \geq 1$, and each
component of $\partial S$ contains at least one of the $v_i$.

\begin{definition}
An ideal triangulation of $S$ is a triangulation of the closed
surface \={S} whose vertex set is exactly \{$v_1,...,v_p$\}.
\end{definition}

Such an ideal triangulation exists if and only if S is not one of
the following surfaces: the sphere with at most two points
removed, the disk with one point removed and the disk with two
points on the boundary removed. We will always assume that S is
not one of these surfaces to insure the existence of an ideal
triangulation.

If $p_{int}$ of the points $v_i$ are in the interior of $\bar{S}$
and if $p_{\partial}$ of the points $v_i$ are on the the boundary $\partial \bar{S}$,
an easy Euler characteristic argument shows that any ideal
triangulation has $n=6g-6+3p_{int} + 2p_{\partial}$ edges.

Two ideal triangulations are considered the same if they are
isotopic. In addition, we require that each ideal triangulation
$\lambda$ is endowed with an indexing of its edges
$\lambda_1,...,\lambda_n$. Let $\Lambda(S)$ denote the set of
isotopy classes of such indexed ideal triangulations $\lambda$.

The set $\Lambda(S)$ admits a natural action of the group $S_n$ of
permutations of $n$ elements, acting by permuting the indices of
the edges of $\lambda$. Namely $\hat{\lambda} = \sigma(\lambda)$
for $\sigma \in S_n$, if its $i$-th edge $\hat{\lambda}_i$ is
equal to $\lambda_{\sigma(i)}$.

Another important transformation of $\Lambda(S)$ is provided by
the $i-th$ $ diagonal$ $ exchange$ $ map$ $\Delta_i : \Lambda(S)
\rightarrow \Lambda(S)$ defined as follows. The $i$-th edge
$\lambda_i$ of an ideal triangulation $\lambda \in \Lambda(S)$ is
adjacent to two triangles. If these two triangles are distinct,
their union forms a square $Q$ with diagonal $\lambda_i$. Then
$\Delta_i(\lambda)$ is obtained from $\lambda$ by replacing edge
$\lambda_i$ by the other diagonal $\hat{\lambda}_i$ of the square
$Q$. By convention, $\Delta_i(\lambda) = \lambda$ when the two
sides of $\lambda_i$ belong to the same triangle; this happens
exactly when $\lambda_i$ is the only edge of $\lambda$ leading to
a puncture of S.

\section{The Exponential Shear Coordinates For the Enhanced Teichm\"uller
Space}

\begin{definition}
The Teichm\"{u}ller space of $S$ is the space $\mathcal{T}(S)$ of
complete hyperbolic metrics on $S$ for which $\partial S$ is
geodesic, considered up to isotopy.
\end{definition}

Consider a complete hyperbolic metric $m \in \mathcal{T}(S)$. It
is well-known that the ends of the complete hyperbolic surface
$(S,m)$ can be of three types: spikes bounded on each side by two
components of $\partial S$ (possibly equal), finite area cusps
bounded on one side by a horocycle; and infinite area funnels
bounded on one side by a simple closed geodesic.

It is convenient to enhance the hyperbolic metric $m \in
\mathcal{T}(S)$ with some additional data, consisting of an
orientation for each closed geodesic bounding a funnel end. Let
the  \textit{enhanced Teichm\"{u}ller space}
$\tilde{\mathcal{T}}(S)$ consist of all isotopy classes of
hyperbolic metrics $m \in \mathcal{T}(S)$ enhanced with such a
choice of orientation. The enhanced Teichm\"{u}ller space
$\tilde{\mathcal{T}}(S)$ inherits from the topology of
$\mathcal{T}(S)$ a topology for which the natural projection
$\tilde{\mathcal{T}}(S) \rightarrow \mathcal{T}(S)$ is a branched
covering map.

Thurston associated a certain
system of coordinates for the enhanced Teichm\"{u}ller space
$\tilde{\mathcal{T}}(S)$ to an ideal triangulation $\lambda$ , called the shear coordinates.

Consider an enhanced hyperbolic metric $m \in
\tilde{\mathcal{T}}(S)$ together with an ideal triangulation
$\lambda$. Each edge $\lambda_i$ specifies a proper homotopy class
of paths going from one end of $(S,m)$ to another end. This proper
homotopy class is also realized by a unique $m$-geodesic $g_i$
such that each end of $g_i$, either converges toward a spike, or
converges towards a cusp end of $S$, or spirals around a closed
geodesic bounding a funnel end in the direction specified by the
enhancement of $m$. The closure of the union of the $g_i$ forms an
$m$-geodesic lamination $g$.

The enhanced hyperbolic metric $m \in \tilde{\mathcal{T}}(S)$ now
associates to an edge $\lambda_i$ of $\lambda$ a positive number
$x_i$ defined as follows. The geodesic $g_i$ separates two
triangle components $T^1_i$ and $T^2_i$ of $S - g$. Isometrically
identify the universal covering of $(S,m)$ to the upper half-space
model $\mathbb{H}^2$ for the hyperbolic plane. Lift $g_i$, $T_i^1$
and $T_i^2$ to a geodesic $\tilde{g_i}$ and the two triangles
$\tilde{\mathcal{T}}^1_i$ and $\tilde{\mathcal{T}}^2_i$ in
$\mathbb{H}^2$ so that the union $\tilde{g_i} \cup
\tilde{\mathcal{T}}^1_i \cup \tilde{\mathcal{T}}^2_i$ forms a
square $\tilde{Q}$ in $\mathbb{H}^2$. Let $z_- , z_+ , z_r , z_l$
be the vertices of $\tilde{Q}$, indexed in such a way that
$\tilde{g_i}$ goes from $z_-$ to $z_+$ and that, for this
orientation of $\tilde{g_i}$ , the points $z_r , z_l$ are
respectively to the right and to the left of $\tilde{g_i}$ for the
orientation of $\tilde{Q}$ given by the orientation of $S$. Then
\begin{equation}
x_i = -\mathrm{crossratio}(z_r , z_l , z_- , z_+) = -\frac{(z_r -
z_-)(z_l - z_+)}{(z_r - z_+)( z_l - z_-)} \nonumber
\end{equation}
Note that $x_i$ is positive since the points $z_l , z_- , z_r ,
z_+$ occur in this order in the real line bounding the upper
half-space $\mathbb{H}^2$.

The real numbers $x_i$ are the \emph{exponential shear
coordinates} of the enhanced hyperbolic metric $m \in
\tilde{\mathcal{T}}(S)$ \cite{Thur1}. The standard shear
coordinates are their logarithms,  $ log(x_i)$,  but the $x_i$
turn out to be better behaved for our purposes.

There is an inverse construction which associates a hyperbolic
metric to each system of positive weights $x_i$ attached to the
edges $\lambda_i$ of the ideal triangulation $\lambda$: Identify
each of the components of $S - \lambda$ to a triangle with
vertices at infinity in $\mathbb{H}^2$ , and glue these hyperbolic
triangles together in such a way that adjacent triangles for a
square whose vertices have cross-ratio $x_i$ as above. This
defines a possibly incomplete hyperbolic metric on the surface
$S$. An analysis of this metric near the ends of $S$ shows that
its completion is a hyperbolic surface $S'$ with a geodesic
boundary, and that each end of an edge of $\lambda$ either spirals
towards a component of $\partial {S'}$ or converges towards a cusp
end of $S'$. Extending $S'$ to a complete hyperbolic metric $m$ on
$S$ whose convex core is isometric to $S'$. In addition, the
spiraling pattern of the ends of $\lambda$ provides an enhancement
of the hyperbolic metric $m$.

The $x_i$ then defines a homeomorphism $\phi_\lambda :
\tilde{\mathcal{T}}(S) \rightarrow \mathbb{R}_+^n$ between
enhanced Teichm\"{u}ller space $\tilde{\mathcal{T}}(S)$ and
$\mathbb{R}_+^n$.

\section{Trace Functions}

A simple closed curve $\alpha$ on the surface $S$ determines a
trace function $T_{\alpha} : \mathcal{T}(S) \rightarrow
\mathbb{R}$, defined as follows: The monodromy of $m \in
\mathcal{T}(S)$ is a group homomorphism $r_m: \pi_1 (S)
\rightarrow PSL_2(\mathbb{R})$ well defined up to conjugation. The
trace of $r_m(\alpha) \in PSL_2(\mathbb{R})$ is only defined up to
sign. Let $T_{\alpha} (m)=|Tr(r_m(\alpha))|$.

This trace function $T_{\alpha}$ has a nice expression in terms of
shear coordinates. Fix an ideal triangulation $\lambda$, and
consider the associated parametrization
$\phi_{\lambda}:\tilde{\mathcal{T}}(S) \rightarrow
\mathbb{R}_{+}^n$ by shear coordinates.

\begin{proposition}
 For every ideal triangulation
$\lambda$ and every simple closed curve $\alpha$, the function
$T_{\alpha} \circ \phi^{-1}_{\lambda}:\mathbb{R}_{+}^n \rightarrow
\mathbb{R}$ is a Laurent polynomial in the square roots $\{ x_1^{
\frac{1}{2}}, x_2^{\frac{1}{2}},...,x_n^{\frac{1}{2}} \}$ of the
shear coordinates.
\end{proposition}

\begin{proof}

     As in \cite{CF1} let us introduce the
    ``left'' and ``right'' turn matrices
    $L \equiv \left(\begin{array}{ccc}1&1\\0&1\end{array}\right)$ and
    $R \equiv \left(\begin{array}{ccc}1&0\\1&1\end{array}\right).$
To each edge $\lambda_i$ in $\lambda$ we associate $S(x_i) =
\left(\begin{array}{ccc}x_{i}^{\frac{1}{2}}&0\\0&x_i^{-\frac{1}{2}}\end{array}\right)$
where the coordinate shear along $\lambda_i$ is $x_i$. For a
closed curve $\alpha$ in $S$, choose any point $\varsigma$ on
$\alpha$ and trace once around $\alpha$ until you return to the
point $\varsigma$. Looking at the directed path traced along
$\alpha$, we record an $S(x_i)$ every time $\alpha$ crosses
$\lambda_i$. If the directed path traced along $\alpha$ crosses
$\lambda_i$ and then $\lambda_j$, we record an $L$ if both
$\lambda_i$ and $\lambda_j$ are asymptotic to each other on the
left of the directed path, and we record a $R$ if both $\lambda_i$
and $\lambda_j$ are asymptotic to each other on the right of the
directed path. This yields a string of matrices
$P_1S(x_{i_1})P_2S(x_{i_2})...P_nS(x_{i_n})$ where the $P_i$ are
either $R$ or $L$ depending on the criterion above.

An argument in \cite{CF1} shows that, up to conjugation, $r_m(\alpha) =
P_1S(x_{i_1})P_2S(x_{i_2})...P_nS(x_{i_n})$.

Note that the trace of
$P_1S(x_{i_1})P_2S(x_{i_2})...P_nS(x_{i_n})$ is a Laurent
poynomial in the $x_i^{\frac{1}{2}}$ with positive coefficients.
In particular, it is positive. Therefore, $T_{\alpha}(m)=
Tr(P_1S(x_{i_1})P_2S(x_{i_2})...P_nS(x_{i_n}))$ is a Laurent
polynomial in the $x_i^{\frac{1}{2}}$.

 \end{proof}

\section{The Chekhov-Fock Algebra}

We consider a quantization of the enhanced Teichm\"{u}ller space
$\tilde{\mathcal{T}}(S)$, by defining a deformation depending on a
parameter $q$, of the algebra $\mathrm{Rat}\tilde{\mathcal{T}}(S)$ of
all the rational functions of $\tilde{\mathcal{T}}(S)$.

Fix an indexed ideal triangulation $\lambda \in \Lambda(S)$. Its
complement $S-\lambda$ has $2n$ spikes converging towards the
punctures, and each spike is delimited by one of the indexed edges
$\lambda_{i}$ of $\lambda$ on one side, and one $\lambda_{j}$ on
the other side; here $i=j$ is possible. For $i,j\in {1,2,...,n}$,
let $\alpha_{ij}^{\lambda}$ denote the number of spikes of
$S-\lambda$ which are delimited on the left by $\lambda_{i}$ and
on the right by $\lambda_{j}$, and set
\begin{equation}
\sigma^{\lambda}_{ij}=\alpha_{ij}^{\lambda}-\alpha_{ji}^{\lambda}
\nonumber.
\end{equation}
Notice that $\sigma^{\lambda}_{ij}$ can only belong to the set
$\{-2,-1,0,1,2\}$, and that
$\sigma^{\lambda}_{ij}=-\sigma^{\lambda}_{ji}$.

 The \textit{Chekhov-Fock algebra} associated to the ideal
triangulation $\lambda$ is the algebra $\mathcal{T}_{\lambda}^{q}$
defined by the generators $X_{1}, X_{1}^{-1}, X_{2}, X_{2},
X_{2}^{-1},...,X_{n},X_{n}^{-1}$, with each pair $X_{i}^{\pm 1}$
associated to an edge $\lambda_{i}$ of $\lambda$, and by the
relations
\begin{equation*} \tag{1.1}
X_{i}X_{j}=q^{2\sigma_{ij}^{\lambda}}X_{j}X_{i}
\end{equation*}
\begin{equation*} \tag{1.2}
X_{i}X_{i}^{-1}=X_{i}^{-1}X_{i}=1
\end{equation*}

If $q=1$, then the $X_i$ commutes and is equal to the Thurston shear
coordinates $x_i$ introduced in Section 2.

The Chekhov-Fock algebra is a Noetherian ring and a right Ore domain so we can
introduce the fraction division algebra
$\hat{\mathcal{T}}_{\lambda}^{q}$ consisting of all formal
fractions $PQ^{-1}$ with $P,Q \in  \mathcal{T}_{\lambda}^{q}$ and
$Q \neq 0$. Two such fractions $P_1Q_1^{-1}$ and $P_2Q_2^{-1}$ are
identified if there exists $S_1,S_2 \in \mathcal{T}_{\lambda}^{q}$
such that $P_1S_1=S_2P_2$ and $Q_1S_1=S_2Q_2$.

Chekhov-Fock \cite{CF1} and Kashaev \cite{Kash1} (see also Liu
\cite{Liu1}) introduced the following:

\begin{theorem}
(Chekhov-Fock, Kashaev, Liu) There exists a family of algebra
isomorphisms

\begin{equation}
\Phi_{\lambda \lambda'}:\hat{\mathcal{T}}_{\lambda'}^{q}
\rightarrow \hat{\mathcal{T}}_{\lambda}^{q} \nonumber
\end{equation}
indexed by pairs of ideal triangulations $\lambda$, $\lambda'$
$\in \Lambda(S)$, which satisfy the following conditions:
\begin{enumerate}
\item $\Phi_{\lambda \lambda''} = \Phi_{\lambda \lambda'} \circ
\Phi_{\lambda' \lambda''}$ for any $\lambda$, $\lambda'$, and
$\lambda''$ $\in \Lambda(S)$.

\item If $\lambda' = \sigma \lambda$ is obtained by reindexing
$\lambda$ by a permutation $\sigma \in S_n$, then $\Phi_{\lambda
\lambda'}(X_{i}^{'})=X_{\sigma(i)}^{'}$ for any $1 \leq i \leq n$.

\end{enumerate}
\end{theorem}

The $\Phi_{\lambda \lambda'}$ are called the Chekhov-Fock
coordinate change isomorphisms. We can now define the quantum
Teichm\"{u}ller space by using the Chekhov-Fock fraction algebras
$\hat{\mathcal{T}}_{\lambda}^{q}$ as charts and the Chekhov-Fock
isomorphisms as coordinate change maps. More precisely:

\begin{definition}
The quantum (enhanced) Teichm\"{u}ller space of a surface S is the
algebra
\begin{equation*}
\mathcal{T}_S^q = \left( \bigsqcup_{\lambda \in \Lambda(S)}
\mathcal{T}_{\lambda}^q \right) / \thicksim \qquad ,
\end{equation*}
where the relation $\thicksim$ is defined by the property that,
for $X \in \mathcal{T}_{\lambda}^q$ and $X' \in
\mathcal{T}_{\lambda '}^q$,
\begin{equation*}
X \thicksim X' \Leftrightarrow X= \Phi_{\lambda \lambda'}(X')
\end{equation*}

\end{definition}

The construction is specially designed so that, when $q=1$, there
is a natural isomorphism between the $\mathcal{T}_S^1$ and the algebra
$\mathrm{Rat} \, \tilde{\mathcal{T}}(S)$ of rational functions on
the enhanced Teichm\"{u}ller space $\tilde{\mathcal{T}}(S)$. See
\cite{CF1}, \cite{BonLiu1}, \cite{Liu1}.

\section{Square Roots}

In the non-quantum case the formula which defines the traces
involves square roots of the shear coordinates. Therefore we need an
algebra which is generated by the square roots of the generators of
the Chekhov-Fock algebra. This leads us to the square root algebra
$\mathcal{T}_{\lambda}^{q^{\frac{1}{4}}}$ defined by the
generators $Z_{1}, Z_{1}^{-1}, Z_{2},
Z_{2}^{-1},...,Z_{n},Z_{n}^{-1}$, where $Z_i = X_i^{\frac{1}{2}}$,
with each pair $Z_{i}^{\pm 1}$ associated to an edge $\lambda_{i}$
of $\lambda$, and by the relations

\begin{equation*} \label{SkewRel1} \tag{1.1}
Z_{i}Z_{j}=q^{\frac{1}{2}\sigma_{ij}^{\lambda}}Z_{j}Z_{i}
\end{equation*}
\begin{equation*} \tag{1.2}
Z_{i}Z_{i}^{-1}=Z_{i}^{-1}Z_{i}=1
\end{equation*}

The square root algebra $\mathcal{T}_{\lambda}^{q^{\frac{1}{4}}}$
is just the Checkhov-Fock algebra with a different $q$. In
particular, we need to choose a $4$-th root, $q^{\frac{1}{4}}$, for
$q$. There is a natural inclusion map:
$$ i:\mathcal{T}_{\lambda}^{q} \hookrightarrow \mathcal{T}_{\lambda}^{q^{\frac{1}{4}}}$$
$$ X_i \mapsto Z_i^2$$
which induces the inclusion:
$$\hat{i}: \hat{\mathcal{T}}_{\lambda}^{q} \hookrightarrow \hat{\mathcal{T}}_{\lambda}^{q^{\frac{1}{4}}}$$
of the fraction division algebras
$\hat{\mathcal{T}}_{\lambda}^{q}$ and
$\hat{\mathcal{T}}_{\lambda}^{q^{\frac{1}{4}}}$.

Unfortunately there is no nice extension of the Chekhov-Fock
coordinate changes to the square root algebra
$\mathcal{T}_{\lambda}^{q^{\frac{1}{4}}}$. This leads us to
introduce the following definitions. The first definition is specially
designed to address the problem that we are facing and the second
definition is very classical.

\begin{definition}

For an ideal triangulation $\lambda$ with edges $\lambda_1$,
$\lambda_2$,...,$\lambda_n$ and a simple closed curve $\alpha$
which crosses edges $\lambda_{i_1}, \lambda_{i_2}, ...
,\lambda_{i_h}$, an element $T$ of
$\hat{\mathcal{T}}_{\lambda}^{q^{\frac{1}{4}}}$ is
\emph{$\alpha$-odd} if it can be written as
$$ T= Z_{i_1}^{-1}Z_{i_2}^{-1}...Z_{i_h}^{-1}R $$
with $ R \in \hat{\mathcal{T}}_{\lambda}^{q} $. The set of
$\alpha$-odd elements is denoted by
$\hat{\mathcal{T}}_{\lambda}^{q^{\frac{1}{4}}}(\alpha)$.

\end{definition}

\begin{remark}

It is worth noting that the set
$\hat{\mathcal{T}}_{\lambda}^{q^{\frac{1}{4}}}(\alpha)$ is not an
algebra. Also for every $T \in
\hat{\mathcal{T}}_{\lambda}^{q^{\frac{1}{4}}}(\alpha)$, the square
$T^2$ is in the subalgebra $ \hat{\mathcal{T}}_{\lambda}^{q}
\subset \hat{\mathcal{T}}_{\lambda}^{q^{\frac{1}{4}}}$

\end{remark}

\begin{definition}
For a monomial $Z_{i_1} Z_{i_2} ... Z_{i_r} \in
\hat{\mathcal{T}}_{\lambda}^{q^{\frac{1}{4}}}$ the \emph{Weyl
ordering coefficient} associated to this monomial is the
coefficient $q^w$ with $w=-\frac{1}{4} \sum_{j < k} \sigma_{i_j
i_k}$.
\end{definition}

The exponent $w$ is engineered so that the quantity $q^w
Z_{i_1}Z_{i_2}...Z_{i_r}$ is unchanged when one permutes the
$Z_{i_s}$'s.

Given an ideal triangulation $\lambda$ of a surface $S$ of genus
$g$ with $p$ punctures, $p_{int}$ on the interior and
$p_{\partial}$ on the boundary, label its triangles by $T_1$,
$T_2$,...,$T_m$, where $m=-2 \chi(S) +
p_{\partial}=4g-4+2p_{int}+p_{\partial}$. Each triangle $T_m$
determines a \emph{triangle algebra}
$\mathcal{T}_{T_m}^{q^{\frac{1}{4}}}$, defined by the generators
$Z_{i,m}$, $Z_{i,m}^{-1}$, $Z_{j,m}$, $Z_{j,m}^{-1}$, $Z_{k,m}$,
$Z_{k,m}^{-1}$ with each pair $Z_{l,m}^{\pm 1}$ associated to an
edge $\lambda_{l}$ of the triangle $T_m$, and by the relations

\begin{equation*}
Z_{i_1,m}Z_{i_2,m}=q^{\frac{1}{2}}Z_{i_2,m}Z_{i_1,m}
\end{equation*}

if $Z_{i_1,m}$,$Z_{i_2,m}\in\mathcal{T}_{T_m}^{q^{\frac{1}{4}}}$
are the generators associated to two sides of $T_m$ with
$Z_{i_1,m}$ associated to the side that comes first when going
counterclockwise at their common vertex.

The square root algebra $\mathcal{T}_{\lambda}^{q^{\frac{1}{4}}}$
has a natural embedding into the tensor product algebra
$\bigotimes^m_{i=1} \mathcal{T}_{T_i}^{q^{\frac{1}{4}}} =
\mathcal{T}_{T_1}^{q^{\frac{1}{4}}} \otimes ... \otimes
\mathcal{T}_{T_m}^{q^{\frac{1}{4}}}$ defined as follows. If the
generator $Z_i$ of $\mathcal{T}_{\lambda}^{q^{\frac{1}{4}}}$ is
associated to the $i$-th edge $\lambda_i$ of $\lambda$, define
\begin{enumerate}
\item $Z_i \rightarrow Z_{i,j} \otimes Z_{i,k}$ if $\lambda_i$
separates two distinct faces $T_j$ and $T_k$, and if $Z_{i,j} \in
\mathcal{T}_{T_j}^{q^{\frac{1}{4}}}$ and $Z_{i,k} \in
\mathcal{T}^{q^{\frac{1}{4}}}_{T_k}$ are the generators associated
to the sides of $T_j$ and $T_k$ corresponding to $\lambda_i$.
\item $Z_i \rightarrow
q^{-\frac{1}{4}}Z_{i_1,j}Z_{i_2,j}=q^{\frac{1}{4}}Z_{i_2,j}Z_{i_1,j}$
if $\lambda_i$ corresponds to the two sides of the same face
$T_j$, and if
$Z_{i_1,j}$,$Z_{i_2,j}\in\mathcal{T}_{T_j}^{q^{\frac{1}{4}}}$ are
the generators associated to these two sides with $Z_{i_1,j}$
associated to the side that comes first when going
counterclockwise at their common vertex.
\end{enumerate}

\begin{figure}[htb]
\centering
 \SetLabels
 ( 0.23 * 0.52) \tiny{(a)}\\
  ( .23 * 1.03) \tiny{$\lambda_j$}\\
  (0.26 * 0.8) \tiny{$\lambda_i$}\\
  ( .23 * 0.63) \tiny{$\lambda_l$}\\
  (0.29 * 0.92) \tiny{$\alpha$}\\
  ( 0.57 * 0.52) \tiny{(b)}\\
  ( .91 * 0.52) \tiny{(c)}\\
  ( 0.23 * -0.07) \tiny{(d)}\\
  ( .57 * -0.07) \tiny{(e)}\\
  ( 0.91 * -0.07) \tiny{(f)}\\
  ( 0.19 * 0.34) \tiny{$\alpha$}\\
  ( 0.57 * 0.27) \tiny{$\alpha$}\\
  ( 0.63 * 0.7) \tiny{$\alpha$}\\
  ( 0.93 * 0.27) \tiny{$\alpha$}\\
  ( 0.84 * 0.67) \tiny{$\alpha$}\\
  ( .36 * 0.77) \tiny{$\lambda_k$}\\
  ( 0.1 * 0.77) \tiny{$\lambda_m$}\\
   ( .61 * 0.86) \tiny{$T_1$}\\
  ( .52 * 0.7) \tiny{$T_2$}\\
\endSetLabels
\centerline{\AffixLabels{
\includegraphics{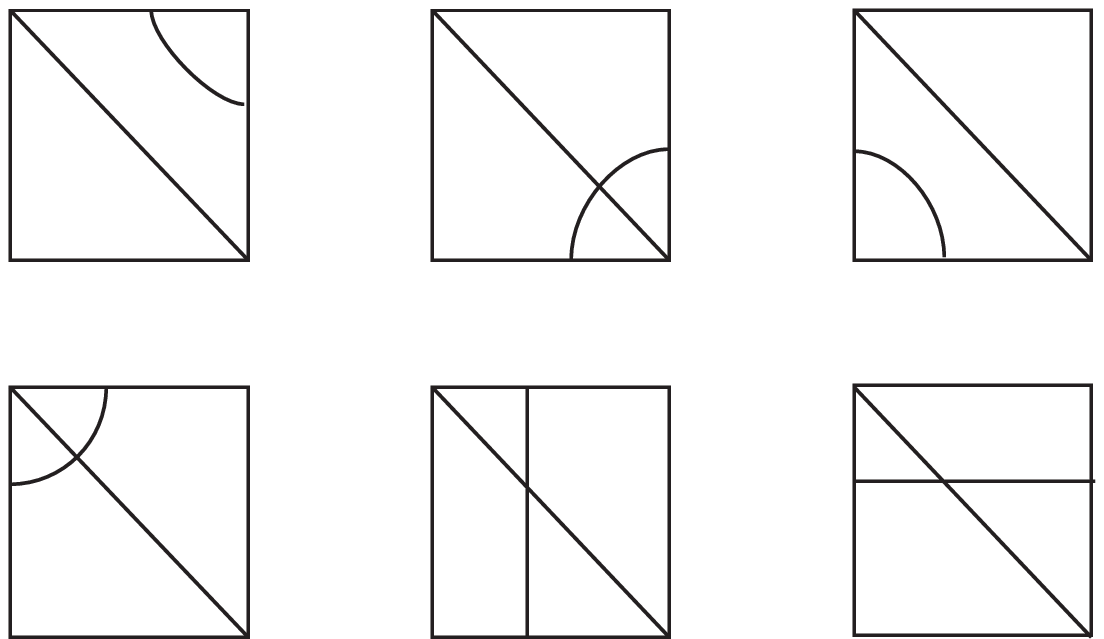}}}
\caption{Cases For Definition of Coordinate Change Maps}
\label{fig:SmallBoxes}
\end{figure}


Consider an ideal triangulation $\lambda$ with edges $\lambda_1 $,
$\lambda_2$,...,$\lambda_n$, and let $\alpha$ be a simple closed
curve in $S$ which is transverse to $\lambda$, where $\alpha$ does not backtrack over the edges of $\lambda$. Namely, $\alpha$ never enters a triangle of
$S-\lambda$ through one side and exits through the same side.

Now consider a square $Q$ in the triangulation $\lambda$ formed by
the edges $\lambda_i$, $\lambda_j$, $\lambda_k$, $\lambda_l,$ and
$\lambda_m$ as in Figure ~\ref{fig:SmallBoxes} (a). Then $\alpha$
can cross the square $Q$ several times. There are six
possibilities for doing so, which are depicted in Figure
~\ref{fig:SmallBoxes}. To each time $\alpha$ crosses $Q$, we
associate a "block" $B \in \bigotimes_{i=1}^m
\mathcal{T}_{T_i}^{q^{\frac{1}{4}}}$ defined as follows.

\begin{enumerate}

\item $B=q^{-\frac{1}{4}} Z_{j,1}^{-1} Z_{k,1}^{-1} $ when
$\alpha$ crosses $Q$ is as in
    Figure ~\ref{fig:SmallBoxes} (a).

\item $B= q^{-\frac{1}{2}} Z_{k,1}^{-1} Z_i^{-1}
\hat{Z}_{l,2}^{-1}$ when $\alpha$ crosses $Q$ is as in
    Figure ~\ref{fig:SmallBoxes} (b).

\item $B= q^{-\frac{1}{4}} Z_{l,2}^{-1} Z_{m,2}^{-1} $ when
$\alpha$ crosses $Q$ is as in
    Figure ~\ref{fig:SmallBoxes} (c).

\item $ B= q^{\frac{1}{2}} Z_{j,1}^{-1} Z_i^{-1} Z_{m,2}^{-1}$
when $\alpha$ crosses $Q$ is as in
    Figure ~\ref{fig:SmallBoxes} (d).

\item $B= Z_{j,1}^{-1} Z_i^{-1} Z_{l,2}^{-1}$ when $\alpha$
crosses $Q$ is as in
    Figure ~\ref{fig:SmallBoxes} (e).

\item $B= Z_{k,1}^{-1} Z_i^{-1} Z_{m,2}^{-1}$ when $\alpha$
crosses $Q$ is as in
    Figure ~\ref{fig:SmallBoxes} (f).

\end{enumerate}

Note the Weyl ordering of the $q$ coefficent of the block $B$ in
each case.

\begin{lemma}
\label{lem:BlockDecomp} Let $B_1,B_2,...,B_r$ be blocks associated
to a simple closed curve $\alpha$ in $S$ crossing squares $Q_1,
Q_2,...Q_r$ as in Figure ~\ref{fig:SmallBoxes}.


 Then every $T \in
\hat{\mathcal{T}}_{\lambda}^{q^{\frac{1}{4}}}(\alpha)$ can be
written in a unique way as $T=q^w B_1B_2...B_rR$ with $R \in
\hat{\mathcal{T}}_{\lambda}^q$, where $q^w$ is the Weyl ordering
coefficient of the blocks $B_i$.

\end{lemma}

\begin{proof}
Every element $T \in
\hat{\mathcal{T}}_{\lambda}^{q^{\frac{1}{4}}}(\alpha)$ can be
written as $ T= Z_{i_1}^{-1}Z_{i_2}^{-1}...Z_{i_h}^{-1}R $ with $
R \in \hat{\mathcal{T}}_{\lambda}^{q} $. The result follows from
the fact that
$B_1B_2...B_r=q^bZ_{i_1}^{-1}Z_{i_2}^{-1}...Z_{i_h}^{-1}$ for some
power of $q$.

\end{proof}


We now want to generalize to the non-commutative context the
coordinate change isomorphisms $\Phi_{\lambda \lambda'}^q$ from
\cite{Liu1}, described in the previous section, by introducing
appropriate linear maps $\Theta^{q}_{\hat{\lambda} \lambda} :
\hat{\mathcal{T}}_{\lambda}^{q^{\frac{1}{4}}}(\alpha) \rightarrow
\hat{\mathcal{T}}_{\hat{\lambda}}^{q^{\frac{1}{4}}}(\alpha)$ in
the following way:

\begin{definition}
\label{Def:LinearMap}  Given a simple closed curve $\alpha$ in the
surface $S$, and $\alpha$-odd ideal triangulations $\lambda$,
$\hat{\lambda}$ separated by a single diagonal exchange, define
$\Theta^{q}_{\hat{\lambda}
\lambda}:\hat{\mathcal{T}}_{\lambda}^{q^{\frac{1}{4}}}(\alpha)
\rightarrow
\hat{\mathcal{T}}_{\hat{\lambda}}^{q^{\frac{1}{4}}}(\alpha)$ as
follows:

If $T=q^wB_1B_2...B_rR$ as in Lemma ~\ref{lem:BlockDecomp},
$\Theta^{q}_{\hat{\lambda} \lambda}(T)$ is obtained from $T$ by

\begin{enumerate}
\item Keeping the same coefficient $q^w$

\item Replacing $R$ with $\Phi_{\hat{\lambda} \lambda}(R) \in
\hat{\mathcal{T}}_{\hat{\lambda}}^q$

\item Replacing $B_i=q^{-\frac{1}{4}} Z_{j,1}^{-1} Z_{k,1}^{-1} $
with $\hat{B}_i=q^{-\frac{1}{2}} \hat{Z}_{j,1}^{-1} \hat{Z}_i^{-1}
\hat{Z}_{k,2}^{-1}$ when the block $B_i$\\
    is associated to the configuration of Figure ~\ref{fig:SmallBoxes} (a).

\item Replacing $B_i= q^{-\frac{1}{2}} Z_{k,1}^{-1} Z_i^{-1}
Z_{l,2}^{-1}$ with
    $\hat{B}_i=q^{-\frac{1}{4}} \hat{Z}_{k,2}^{-1} \hat{Z}_{l,2}^{-1} $ when the block
    $B_i$\\
    is associated to the configuration of
    Figure ~\ref{fig:SmallBoxes} (b).

\item Replacing $B_i=q^{-\frac{1}{4}} Z_{l,2}^{-1} Z_{m,2}^{-1} $
with
    $\hat{B}_i=q^{-\frac{1}{2}} \hat{Z}_{l,2}^{-1} \hat{Z}_i^{-1}  \hat{Z}_{m,1}^{-1}$ when the block
    $B_i$\\
    is associated to the configuration of
    Figure ~\ref{fig:SmallBoxes} (c).

\item Replacing  $ B_i=q^{\frac{1}{2}}Z_{j,1}^{-1} Z_i^{-1}
Z_{m,2}^{-1}$ with
    $\hat{B}_i=q^{\frac{1}{4}}  \hat{Z}_{j,1}^{-1} \hat{Z}_{m,1}^{-1}$
    when the block $B_i$\\
    is associated to the configuration of
    Figure ~\ref{fig:SmallBoxes} (d).

\item Replacing $B_i=Z_{j,1}^{-1} Z_i^{-1} Z_{l,2}^{-1}$ with
$\hat{B}_i=\hat{Z}_{j,1}^{-1}  (\hat{Z}_i+ \hat{Z}_i^{-1})^{-1}
\hat{Z}_{l,2}^{-1}$ when the block $B_i$\\
    is associated to the configuration of
    Figure ~\ref{fig:SmallBoxes} (e).

\item Replacing $B_i=Z_{k,1}^{-1}  Z_i^{-1} Z_{m,2}^{-1}$ with
    $ \hat{B}_i=\hat{Z}_{k,2}^{-1} (\hat{Z}_i+ \hat{Z}_i^{-1}) \hat{Z}_{m,1}^{-1} $
when the block $B_i$\\
    is associated to the configuration of
    Figure ~\ref{fig:SmallBoxes} (f).

\end{enumerate}

\end{definition}

Remark: $\Theta_{\hat{\lambda} \lambda}^q$ is only a linear map,
not an algebra homomorphism. Indeed,
$\hat{\mathcal{T}}_{\lambda}^{q^{\frac{1}{4}}}(\alpha)$ is not
even an algebra.

\begin{lemma}
\label{Lem:StayOdd}

$\Theta^{q}_{\hat{\lambda} \lambda}(T)$ is $\alpha$-odd.

\end{lemma}
\begin{proof}
The only case which requires some thought is that of blocks of
type (7) and (8) from Definition ~\ref{Def:LinearMap}. However,
note that in type (7)
\begin{equation*}\hat{B}_i=\hat{Z}_{j,1}^{-1} (\hat{Z}_i+
\hat{Z}_i^{-1})^{-1} \hat{Z}_{l,2}^{-1}=\hat{Z}_{j,1}^{-1}
\hat{Z}_i^{-1} (1+ \hat{Z}_i^{-2})^{-1}
\hat{Z}_{l,2}^{-1}=\hat{Z}_{j,1}^{-1} \hat{Z}_i^{-1}
\hat{Z}_{l,2}^{-1} (1+ q^{-1}\hat{Z}_i^{-2})^{-1} \end{equation*}
with type (8) working the same way.
\end{proof}

\begin{lemma}
\label{Lem:SkewComRel}

The map $\Theta^{q}_{\hat{\lambda} \lambda} :
\hat{\mathcal{T}}_{\lambda}^{q^{\frac{1}{4}}}(\alpha) \rightarrow
\hat{\mathcal{T}}_{\hat{\lambda}}^{q^{\frac{1}{4}}}(\alpha)$ is
independent of the order of the blocks $B_i$.

\end{lemma}

\begin{proof}
Note that, when two blocks $B_1$, $B_2$ are replaced by blocks
$\hat{B}_1$ and $\hat{B}_2$, then $\hat{B}_1$ and $\hat{B}_2$
satisfy the same skew commutativity relation as $B_1$ and $B_2$.
Namely, if $B_1B_2 = q^{2 b}B_2B_1$, then $\hat{B}_1\hat{B}_2 =
q^{2b}\hat{B}_2\hat{B}_1$.

This follows from a  simple computation. For instance,  if $B_1$ and $B_2$
are respectively of type (7) and (8) of
Definition~\ref{Def:LinearMap} then:

\begin{equation*}
\begin{split}
B_1 B_2 &= ( Z_{j,1}^{-1} Z_i^{-1} Z_{l,2}^{-1} ) ( Z_{k,1}^{-1}
Z_i^{-1} Z_{m,2}^{-1} ) \\ & = q^{-1}( Z_{k,1}^{-1} Z_i^{-1}
Z_{m,2}^{-1} ) ( Z_{j,1}^{-1} Z_i^{-1} Z_{l,2}^{-1} ) = q^{-1} B_2
B_1
\end{split}
\end{equation*}
 and
 \begin{equation*}
 \begin{split}
 \hat{B}_1 \hat{B}_2 &=( \hat{Z}_{j,1}^{-1}  (\hat{Z}_i+ \hat{Z}_i^{-1})^{-1} \hat{Z}_{l,2}^{-1} )(
\hat{Z}_{k,2}^{-1} (\hat{Z}_i+ \hat{Z}_i^{-1}) \hat{Z}_{m,1}^{-1}
)
\\ & =
( \hat{Z}_{j,1}^{-1}  \hat{Z}_{l,2}^{-1} )( \hat{Z}_{k,2}^{-1} (\hat{Z}_i+ \hat{Z}_i^{-1})^{-1}  (\hat{Z}_i+ \hat{Z}_i^{-1}) \hat{Z}_{m,1}^{-1} ) \\
& =( \hat{Z}_{j,1}^{-1}  \hat{Z}_{l,2}^{-1} )( \hat{Z}_{k,2}^{-1}
\hat{Z}_{m,1}^{-1} )  = q^{-1} ( \hat{Z}_{k,2}^{-1}
\hat{Z}_{m,1}^{-1} )( \hat{Z}_{j,1}^{-1}  \hat{Z}_{l,2}^{-1} ) \\
& = q^{-1} ( \hat{Z}_{k,2}^{-1} \hat{Z}_{m,1}^{-1} )(
\hat{Z}_{j,1}^{-1} (\hat{Z}_i+ \hat{Z}_i^{-1})(\hat{Z}_i+
\hat{Z}_i^{-1})^{-1} \hat{Z}_{l,2}^{-1} )
\\ & = q^{-1} ( \hat{Z}_{k,2}^{-1} (\hat{Z}_i+ \hat{Z}_i^{-1}) \hat{Z}_{m,1}^{-1} )( \hat{Z}_{j,1}^{-1} (\hat{Z}_i+ \hat{Z}_i^{-1})^{-1} \hat{Z}_{l,2}^{-1} )=q^{-1} \hat{B}_2 \hat{B}_1.
\end{split}
\end{equation*}

The result immediately follows from this property.
\end{proof}

The map $\Theta^q_{\hat{\lambda} \lambda}$ is specially designed
so that:

\begin{lemma}
\label{Lem:SameAsOld}

For all $Z_{i_1}^{- 1} Z_{i_2}^{- 1}...Z_{i_l}^{- 1} R \in
\hat{\mathcal{T}}_{\lambda}^{q^{\frac{1}{4}}}(\alpha)$,

$[\Theta^q_{\hat{\lambda} \lambda}(Z_{i_1}^{- 1} Z_{i_2}^{-
1}...Z_{i_l}^{- 1} R)]^2= \Phi_{\hat{\lambda} \lambda}((Z_{i_1}^{-
1} Z_{i_2}^{- 1}...Z_{i_l}^{- 1} R)^2)$, where $(Z_{i_1}^{- 1}
Z_{i_2}^{- 1}...Z_{i_l}^{- 1} R)^2 \in
\hat{\mathcal{T}}_{\lambda}^q$.

\end{lemma}

\begin{proof}

This lemma follows from simple calculations. Given a simple closed
curve $\alpha$ which crosses edges
$\lambda_{i_1},...,\lambda_{i_l}$, label
$Z_{\lambda}=Z_{i_1}^{-1}...Z_{i_l}^{-1}$. The definition of
$\Theta^q_{\hat{\lambda} \lambda}$ was specifically designed so
that
\begin{equation}
\label{Equ:SmallEquation} \Phi_{\hat{\lambda}
\lambda}((Z_{\lambda})^2) = [\Theta^q_{\hat{\lambda}
\lambda}(Z_{\lambda})]^2.
\end{equation}
For example, consider $Z_{\lambda}= Z_{j,1}^{-1} Z_i^{-1}
Z_{l,2}^{-1}$ as in (7) from Definition ~\ref{Def:LinearMap}. Then
we have:

\begin{equation*}
\begin{split}
[\Theta^q_{\hat{\lambda} \lambda}(Z_{\lambda})]^2 & = [\Theta^q_{\hat{\lambda} \lambda}(Z_{j,1}^{-1} Z_i^{-1} Z_{l,2}^{-1})]^2 \\
& = (\hat{Z}_{j,1}^{-1}  (\hat{Z}_i+ \hat{Z}_i^{-1})^{-1} \hat{Z}_{l,2}^{-1})(\hat{Z}_{j,1}^{-1}  (\hat{Z}_i+ \hat{Z}_i^{-1})^{-1} \hat{Z}_{l,2}^{-1}) \\
& = (\hat{Z}_{j,1}^{-1}  (1 + \hat{Z}_i^{2})^{-1} \hat{Z}_i \hat{Z}_{l,2}^{-1})(\hat{Z}_{j,1}^{-1}  \hat{Z}_i(1 + \hat{Z}_i^{2})^{-1} \hat{Z}_{l,2}^{-1}) \\
& = \hat{Z}_{j,1}^{-2}  (1 + q \hat{Z}_i^{2})^{-1} \hat{Z}_i^2 (1 + q^{-1}\hat{Z}_i^{2})^{-1} \hat{Z}_{l,2}^{-2} \\
& = \hat{Z}_{j,1}^{-2}  (1 + q \hat{Z}_i^{2})^{-1} \hat{Z}_i^2 \hat{Z}_{l,2}^{-2} (1 + q\hat{Z}_i^{2})^{-1}  \\
& = \Phi_{\hat{\lambda} \lambda}(Z_{j,1}^{-2} Z_i^{-2}
Z_{l,2}^{-2}) = \Phi_{\hat{\lambda} \lambda}([Z_{j,1}^{-1}
Z_i^{-1} Z_{l,2}^{-1}]^2)
\end{split}
\end{equation*}

Next we will prove a small lemma:

\begin{sublemma}
\label{Lem:SubLem1} Given a simple closed curve $\alpha$ which
crosses edges $\lambda_{i_1},...,\lambda_{i,n}$ of ideal
triangulation $\lambda$, then $\Phi_{\hat{\lambda}
\lambda}(Z_{\lambda} X_r Z_{\lambda}) = \Theta^q_{\hat{\lambda}
\lambda}(Z_{\lambda}) \Phi_{\hat{\lambda} \lambda}( X_r)
\Theta^q_{\hat{\lambda} \lambda}(Z_{\lambda})$ , for
$Z_{\lambda}=Z_{i_1}^{-1}...Z_{i_l}^{-1}$ and for all $r \in \{
1,2...n \}$.
\end{sublemma}
\begin{proof}
Given that $X_r Z_{\lambda} = q^a Z_{\lambda} X_r$ we have:

\begin{equation*}
\begin{split}
\Phi_{\hat{\lambda} \lambda}(Z_{\lambda} X_r Z_{\lambda}) & =
q^a \Phi_{\hat{\lambda} \lambda}(Z_{\lambda} Z_{\lambda} X_r )  =
q^a \Phi_{\hat{\lambda} \lambda}(Z_{\lambda} Z_{\lambda})\Phi_{\hat{\lambda} \lambda}( X_r ) \\
& = q^a \Theta^q_{\hat{\lambda}
\lambda}(Z_{\lambda})\Theta^q_{\hat{\lambda} \lambda}(
Z_{\lambda})\Phi_{\hat{\lambda} \lambda}( X_r )
=\Theta^q_{\hat{\lambda} \lambda}(Z_{\lambda})\Phi_{\hat{\lambda}
\lambda}( X_r ) \Theta^q_{\hat{\lambda} \lambda}( Z_{\lambda})
\end{split}
\end{equation*} \end{proof}

A direct corollary of Sublemma~\ref{Lem:SubLem1} is,  given a
polynomial $P \in \mathcal{T}_{ \lambda}^q$, $\Phi_{\hat{\lambda}
\lambda}(Z_{\lambda} P Z_{\lambda}) = \Theta^q_{\hat{\lambda}
\lambda}(Z_{\lambda}) \Phi_{\hat{\lambda} \lambda}( P)
\Theta^q_{\hat{\lambda} \lambda}(Z_{\lambda})$. Using this
corollary, we then have, given polynomials $P , Q \in
\mathcal{T}_{ \lambda}^q$:

\begin{equation*}
\begin{split}
\Phi_{\hat{\lambda} \lambda}(Z_{\lambda} P Q^{-1} Z_{\lambda}) & =
\Phi_{\hat{\lambda} \lambda}(Z_{\lambda} P Z_{\lambda})
\Phi_{\hat{\lambda} \lambda}(Z_{\lambda}^{-1} Q^{-1} Z_{\lambda})
  = \Theta^q_{\hat{\lambda} \lambda}(Z_{\lambda})
\Phi_{\hat{\lambda} \lambda}( P) \Theta^q_{\hat{\lambda}
\lambda}(Z_{\lambda}) [\Phi_{\hat{\lambda} \lambda}(Z_{\lambda} Q
Z_{\lambda})]^{-1}
\Phi_{\hat{\lambda} \lambda}(Z_{\lambda}^2) \\
& = \Theta^q_{\hat{\lambda} \lambda}(Z_{\lambda})
\Phi_{\hat{\lambda} \lambda}( P) \Theta^q_{\hat{\lambda}
\lambda}(Z_{\lambda}) [\Theta^q_{\hat{\lambda}
\lambda}(Z_{\lambda})]^{-1} [\Phi_{\hat{\lambda} \lambda}(Q)]^{-1}
   [\Theta^q_{\hat{\lambda} \lambda}(Z_{\lambda})]^{-1}
\Theta^q_{\hat{\lambda} \lambda}(Z_{\lambda})
\Theta^q_{\hat{\lambda} \lambda}(Z_{\lambda}) \text{ }
(\text{Using
}~\ref{Equ:SmallEquation}) \\
&  =\Theta^q_{\hat{\lambda} \lambda}(Z_{\lambda})
\Phi_{\hat{\lambda} \lambda}( P) [\Phi_{\hat{\lambda}
\lambda}(Q)]^{-1} \Theta^q_{\hat{\lambda} \lambda}(Z_{\lambda}) =
\Theta^q_{\hat{\lambda} \lambda}(Z_{\lambda}) \Phi_{\hat{\lambda}
\lambda}( P Q^{-1}) \Theta^q_{\hat{\lambda} \lambda}(Z_{\lambda}).
\end{split}
\end{equation*}

 Now we can finally prove the lemma with the following computations:

\begin{equation*}
\begin{split}
& [\Theta^q_{\hat{\lambda} \lambda}(Z_{i_1}^{- 1} Z_{i_2}^{-
1}...Z_{i_l}^{- 1} R)]^2
 =\Theta^q_{\hat{\lambda}
\lambda}(Z_{i_1}^{- 1} Z_{i_2}^{- 1}...Z_{i_l}^{- 1} R)
 \Theta^q_{\hat{\lambda} \lambda}(Z_{i_1}^{- 1} Z_{i_2}^{- 1}...Z_{i_l}^{- 1} R) \\
 & = \Theta^q_{\hat{\lambda} \lambda}(Z_{i_1}^{- 1} Z_{i_2}^{- 1}...Z_{i_l}^{- 1})\Phi_{\hat{\lambda} \lambda}( R)
\Theta^q_{\hat{\lambda} \lambda}(Z_{i_1}^{- 1} Z_{i_2}^{- 1}...Z_{i_l}^{- 1})\Phi_{\hat{\lambda} \lambda}( R) \\
& =\Phi_{\hat{\lambda} \lambda}(Z_{i_1}^{- 1} Z_{i_2}^{-
1}...Z_{i_l}^{- 1} R
Z_{i_1}^{- 1} Z_{i_2}^{- 1}...Z_{i_l}^{- 1})\Phi_{\hat{\lambda} \lambda}( R) \\
& = \Phi_{\hat{\lambda} \lambda}(Z_{i_1}^{- 1} Z_{i_2}^{-
1}...Z_{i_l}^{- 1} R Z_{i_1}^{- 1} Z_{i_2}^{- 1}...Z_{i_l}^{- 1}
R)
 =\Phi_{\hat{\lambda} \lambda}([Z_{i_1}^{- 1} Z_{i_2}^{-
1}...Z_{i_l}^{- 1} R]^2)
\end{split}
\end{equation*}

This concludes the proof of Lemma~\ref{Lem:SameAsOld}

\end{proof}

This leads us to another lemma.

\begin{lemma}
\label{Lem:InverseMap} If ideal triangulations $\lambda$ and
$\hat{\lambda}$ are separated by a single diagonal exchange then
the maps $\Theta^q_{\lambda \hat{\lambda}}$ and
$\Theta^q_{\hat{\lambda} \lambda }$ are such that $
\Theta^q_{\hat{\lambda} \lambda } = (\Theta^q_{\lambda
\hat{\lambda}})^{-1} $

\end{lemma}

\begin{figure}[htb]
\centering
 \SetLabels
 ( 0.2 * 1.06) $\lambda_j$\\
  ( .8 * 1.06) $\hat{\lambda}_j$\\
  ( .2 * 0.06) $\lambda_l$\\
  ( .8 * .08) $\hat{\lambda}_l$\\
( .16 * 0.56) $\lambda_i$\\
  ( .84 * .56) $\hat{\lambda}_i$\\
  ( -0.04 * 0.5) $\lambda_m$\\
  ( 1.04 * 0.5) $\hat{\lambda}_k$\\
  ( .4 * 0.5) $\lambda_k$\\
  ( .62 * .5) $\hat{\lambda}_m$\\
   ( .11 * 0.75) $T_1$\\
  ( .25 * .3) $T_2$\\
   ( .89 * 0.75) $T_1$\\
  ( .75 * .3) $T_2$\\
\endSetLabels
\centerline{\AffixLabels{
\includegraphics{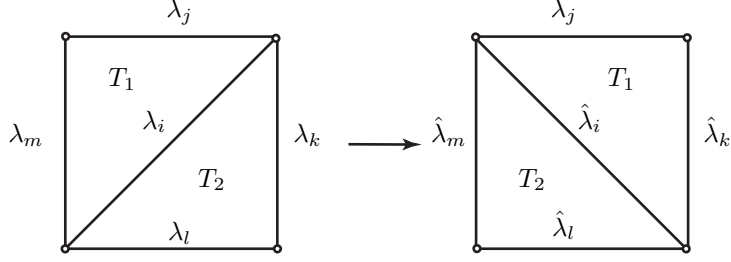}}}
\caption{Diagonal Exchange} \label{Fig:DiagExchangeComputation}
\end{figure}

\begin{proof}
To prove this it is sufficient to show this is true for the six
blocks in Definition~\ref{Def:LinearMap}. Let the edges of
$\lambda$ and $\hat{\lambda}$ involved in the diagonal exchange be
labeled as represented in Figure~\ref{Fig:DiagExchangeComputation}. The
result then follows from computations, all similar to the
following.

\begin{equation*}
\begin{split}
 (\Theta^q_{\lambda \hat{\lambda}})^{-1}(Z_{k,2}^{-1} (Z_i+ Z_i^{-1}) Z_{m,1}^{-1}) & = \hat{Z}_{k,1}^{-1}  \hat{Z}_i^{-1} \hat{Z}_{m,2}^{-1} \\
\Theta^q_{\hat{\lambda} \lambda }(Z_{k,2}^{-1} (Z_i+ Z_i^{-1}) Z_{m,1}^{-1}) & = \Theta^q_{\hat{\lambda} \lambda }(Z_{k,2}^{-1} Z_i^{-1} (1 + Z_i^{2}) Z_{m,1}^{-1}) \\
& = \Theta^q_{\hat{\lambda} \lambda }(Z_{k,2}^{-1} Z_i^{-1}
Z_{m,1}^{-1} (1 + q^{-1} Z_i^{2}))
 = \Theta^q_{\hat{\lambda} \lambda }(Z_{k,2}^{-1} Z_i^{-1}  Z_{m,1}^{-1}) \Phi_{ \hat{\lambda} \lambda}((1 + q^{-1} Z_i^{2})) \\
& = \hat{Z}_{k,1}^{-1} (\hat{Z}_i + \hat{Z}_i^{-1})^{-1}
\hat{Z}_{m,2}^{-1} (1 + q^{-1} \hat{Z}_i^{-2})
 = \hat{Z}_{k,1}^{-1} \hat{Z}_i^{-1}(1 + \hat{Z}_i^{-2})^{-1} \hat{Z}_{m,2}^{-1} (1 + q^{-1} \hat{Z}_i^{-2}) \\
& = \hat{Z}_{k,1}^{-1} \hat{Z}_i^{-1}(1 + \hat{Z}_i^{-2})^{-1} (1
+ \hat{Z}_i^{-2}) \hat{Z}_{m,2}^{-1}
 = \hat{Z}_{k,1}^{-1} \hat{Z}_i^{-1} \hat{Z}_{m,2}^{-1}
\end{split}
\end{equation*}

\end{proof}

The following lemma about the $\Theta^q_{\lambda \hat{\lambda}}$
makes computations easier.

\begin{lemma}
\label{Lem:ComputationForCoordChangeMaps} Given two ideal
triangulations $\lambda$ and $\hat{\lambda}$ which differ only by
a diagonal exchange and if the edges of $\lambda$ and $\hat{\lambda}$ involved in the
diagonal exchange are labeled as depicted in
Figure~\ref{Fig:DiagExchangeComputation}, then the following
relations are satisfied:

\begin{equation*}
\begin{split}
\Theta^q_{\lambda \hat{\lambda}}
(\hat{Z}_{j,1}\hat{Z}_i\hat{Z}_{l,2} +
\hat{Z}_{j,1}^{-1}\hat{Z}_i\hat{Z}_{l,2}^{-1} +
\hat{Z}_{j,1}^{-1}\hat{Z}_i^{-1}\hat{Z}_{l,2}^{-1}) & =
Z_{j,1}Z_iZ_{l,2} + Z_{j,1}Z_i^{-1}Z_{l,2} +
Z_{j,1}^{-1}Z_i^{-1}Z_{l,2}^{-1} \\
 \Theta^q_{\lambda \hat{\lambda}}(\hat{Z}_{k,1}\hat{Z}_i\hat{Z}_{m,2} + \hat{Z}_{k,1}\hat{Z}_i^{-1}\hat{Z}_{m,2}+ \hat{Z}_{k,1}^{-1}\hat{Z}_i^{-1}\hat{Z}_{m,2}^{-1}) & = Z_{k,1}Z_iZ_{m,2}
 + Z_{k,1}^{-1}Z_iZ_{m,2}^{-1} + Z_{k,1}^{-1}Z_i^{-1}Z_{m,2}^{-1} \\
 \Theta^q_{\lambda
\hat{\lambda}}(q^{-\frac{1}{4}}\hat{Z}_{j,1} \hat{Z}_{k,1} +
q^{-\frac{1}{4}}\hat{Z}_{j,1}^{-1} \hat{Z}_{k,1}^{-1}) & =
q^{-\frac{1}{2}}Z_{j,1}Z_iZ_{k,1} + q^{-\frac{1}{2}}Z_{j,1}^{-1}
Z_i^{-1} Z_{k,1}^{-1} \\
 \Theta^q_{\lambda
\hat{\lambda}}(q^{-\frac{1}{2}}\hat{Z}_{k,1} \hat{Z}_{i}
\hat{Z}_{l,2} + q^{-\frac{1}{2}}\hat{Z}_{k,1}^{-1}
\hat{Z}_{i}^{-1} \hat{Z}_{l,2}^{-1}) & =
q^{-\frac{1}{4}}Z_{k,2}Z_{l,2} + q^{-\frac{1}{4}}Z_{k,2}^{-1}
Z_{l,2}^{-1} \\
 \Theta^q_{\lambda
\hat{\lambda}}(q^{\frac{1}{2}}\hat{Z}_{j,1} \hat{Z}_{i}
\hat{Z}_{m,2} + q^{\frac{1}{2}}\hat{Z}_{j,1}^{-1} \hat{Z}_{i}^{-1}
\hat{Z}_{m,2}^{-1}) & = q^{\frac{1}{4}}Z_{j,1}Z_{m,1} +
q^{\frac{1}{4}}Z_{j,1}^{-1} Z_{m,1}^{-1} \\
 \Theta^q_{\lambda
\hat{\lambda}}(q^{-\frac{1}{4}}\hat{Z}_{l,2} \hat{Z}_{m,2} +
q^{-\frac{1}{4}}\hat{Z}_{l,2}^{-1} \hat{Z}_{m,2}^{-1}) & =
q^{-\frac{1}{2}}Z_{l,2}Z_iZ_{m,1} + q^{-\frac{1}{2}}Z_{l,1}^{-1}
Z_i^{-1} Z_{m,2}^{-1} \\
\end{split}
\end{equation*}

\end{lemma}

\begin{proof}
This lemma follows from simple calculations of which we will only
do one.

The first line from the lemma is the result of the following
computations:
\begin{equation*}
\begin{split}
& \Theta^q_{\lambda \hat{\lambda}}(\hat{Z}_{j,1}\hat{Z}_i\hat{Z}_{l,2} + \hat{Z}_{j,1}^{-1}\hat{Z}_i\hat{Z}_{l,2}^{-1}+ \hat{Z}_{j,1}^{-1}\hat{Z}_i^{-1}\hat{Z}_{l,2}^{-1}) \\
& = \Theta^q_{\lambda \hat{\lambda}}(q^{-1}\hat{Z}_{j,1}^2
\hat{Z}_{j,1}^{-1} \hat{Z}_i^{-1} \hat{Z}_{l,2}^{-1}
\hat{Z}_i^2\hat{Z}_{l,2}^2
+ q\hat{Z}_i^2 \hat{Z}_{j,1}^{-1}\hat{Z}_i^{-1}\hat{Z}_{l,2}^{-1}+ \hat{Z}_{j,1}^{-1}\hat{Z}_i^{-1}\hat{Z}_{l,2}^{-1}) \\
& = q^{-1}(1+qZ_i^2)Z_{j,1}^2 Z_{j,1}^{-1} (Z_i + Z_i^{-1})^{-1}Z_{l,2}^{-1}Z_i^{-2}(1+qZ_i^2)Z_{l,2}^2 \\
& \qquad \qquad + qZ_i^{-2}Z_{j,1}^{-1}(Z_i + Z_i^{-1})^{-1} Z_{l,2}^{-1}+ Z_{j,1}^{-1}(Z_i + Z_i^{-1})^{-1} Z_{l,2}^{-1} \\
& = q^{-1}Z_{j,1}(1+qZ_i^2) (1 +
Z_i^2)^{-1}Z_iZ_{l,2}^{-1}(Z_i^{-2}+q)Z_{l,2}^2
 + (1+ qZ_i^{-2})Z_{j,1}^{-1}(Z_i + Z_i^{-1})^{-1} Z_{l,2}^{-1}\\
& =  Z_{j,1}Z_iZ_{l,2} +
q^{-1}Z_{j,1}Z_iZ_{l,2}^{-1}Z_i^{-2}Z_{l,2}^2
+ Z_{j,1}^{-1}(1+ Z_i^{-2})(1 + Z_i^{-2})^{-1} Z_i^{-1}Z_{l,2}^{-1}\\
& = Z_{j,1}Z_iZ_{l,2} + Z_{j,1}Z_i^{-1}Z_{l,2} +
Z_{j,1}^{-1}Z_i^{-1}Z_{l,2}^{-1}
\end{split}
\end{equation*}

\end{proof}

\section{The Pentagon Relation for Square Roots}

The goal of this section is to show that the linear maps
$\Theta^q_{\lambda \hat{\lambda}}$ constructed in the previous
section are compatible with the pentagon relation satisfied by the
diagonal exchange maps, $\Delta_i$, which were introduced in Section 2.

\begin{figure}[htb]
\centering
 \SetLabels
 ( 0.5 * .7) \tiny{$\lambda_a$}\\
 ( 0.21 * .7) \tiny{$\lambda_a$}\\
 ( 0.14 * .08) \tiny{$\lambda_a$}\\
 ( 0.87 * .07) \tiny{$\lambda_a$}\\
  ( .96 * .68) \tiny{$\lambda_b$}\\
  ( .85 * .63) \tiny{$\lambda_x$}\\
  ( .85 * .57) \tiny{$\lambda_y$}\\
  ( 0.83 * .48) \tiny{$\lambda_d$}\\
  ( .74 * .5) \tiny{$\lambda_e$}\\
  ( 1.01 * .5) \tiny{$\lambda_c$}\\
  ( .88 * .40) \tiny{$\lambda^{2}$}\\
  ( .74 * -.04) \tiny{$\lambda^{3}$}\\
  ( 0.28 * -.03) \tiny{$\lambda^{4}$}\\
  ( .14 * .42) \tiny{$\lambda^{0}=\alpha_{x \rightarrow y}\lambda^{5}$}\\
( 0.5 * .65) \tiny{$\lambda^{1}$}\\
( .43 * .88) \tiny{$\alpha$}\\
(0.8 * 0.7) \tiny{$\lambda_a$}\\
  ( 0.51 * .84) \tiny{$\beta$}\\
  ( .69 * .07) \tiny{$\lambda_x$}\\
  ( .72 * .18) \tiny{$\lambda_y$}\\
  ( .21 * .15) \tiny{$\lambda_x$}\\
  ( .27 * .10) \tiny{$\lambda_y$}\\
  ( .145 * .61) \tiny{$\lambda_x$}\\
  ( .08 * .6) \tiny{$\lambda_y$}\\
  ( .45 * .8) \tiny{$\lambda_x$}\\
  ( .5 * .92) \tiny{$\lambda_y$}\\
\endSetLabels
\centerline{\AffixLabels{
\includegraphics{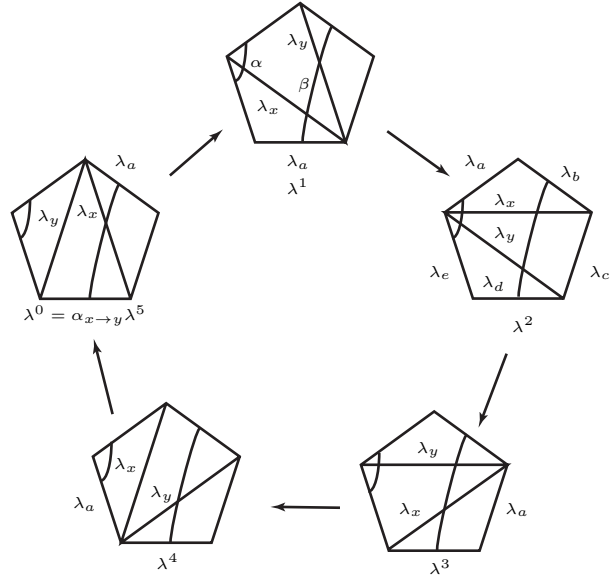}}}
\caption{The Pentagon Relation} \label{fig:PentRelSquares}
\end{figure}

Consider a pentagon cycle of geodesic laminations $\lambda^{0},
\lambda^{1} = \Delta_x(\lambda^{0}), \lambda^{2} =
\Delta_y(\lambda^{1}), \lambda^{3} = \Delta_x(\lambda^{2}),
\lambda^{4} = \Delta_y(\lambda^{3}), \lambda^{5} =
\Delta_x(\lambda^{4}) = \alpha_{x \rightarrow y}(\lambda^{0})$
as represented in Figure ~\ref{fig:PentRelSquares}.

\begin{lemma}
\label{lem:PentagonRelation} The Pentagon Relation
$$\Theta^q_{\lambda^{0} \lambda^{1}} \circ \Theta^q_{\lambda^{1} \lambda^{2}} \circ \Theta^q_{\lambda^{2} \lambda^{3}} \circ
\Theta^q_{\lambda^{3} \lambda^{4}} \circ
\Theta^q_{\lambda^{4} \lambda^{5}} = \mathrm{Id} $$ is
satisfied.

\end{lemma}

\begin{proof}
There are only two non isotopic curves to consider, which are the
$\alpha$ and $\beta$ curves depicted in Figure~\ref{fig:PentRelSquares}.
First we will consider $\alpha$. If we let $\alpha$ be as represented in
Figure~\ref{fig:PentRelSquares} and label the edges of the
pentagons also as depicted in Figure~\ref{fig:PentRelSquares} then we
really only need to look at where $q^{-\frac{1}{4}}Z_{d,1}^{-1}
Z_{e,1}^{-1}$ is mapped to. If we apply the definition of
$\Theta^q_{\lambda \hat{\lambda}}$ and use
Lemma~\ref{Lem:InverseMap} on this monomial we obtain:
\begin{equation*}
\begin{split}
\Theta^q_{\lambda^1 \lambda^0}(q^{-\frac{1}{4}}Z_{d,1}^{-1} Z_{e,1}^{-1}) & = q^{-\frac{1}{2}}Z_{b,2}^{-1} Z_x^{-1} Z_{c,1}^{-1} \\
\Theta^q_{\lambda^2 \lambda^1}(q^{-\frac{1}{2}}Z_{b,2}^{-1}
Z_x^{-1} Z_{c,1}^{-1}) & =
q^{-\frac{1}{4}}Z_{e,3}^{-1}Z_{y,3}^{-1}q^{-\frac{1}{2}}Z_{y,2}^{-1}Z_x^{-1}Z_{a,1}^{-1} \\
& = q^{-\frac{3}{4}}Z_{e,3}^{-1}Z_y^{-1}Z_x^{-1}Z_{a,1}^{-1} \\
\Theta^q_{\lambda^4 \lambda^0}(q^{-\frac{1}{4}}Z_{d,1}^{-1} Z_{e,1}^{-1}) & = q^{-\frac{1}{4}}Z_{a,1}^{-1} Z_{b,1}^{-1} \\
\Theta^q_{\lambda^3 \lambda^4}(q^{-\frac{1}{4}}Z_{a,1}^{-1} Z_{b,1}^{-1}) & = q^{-\frac{1}{2}} Z_{c,2}^{-1} Z_y^{-1} Z_{d,1}^{-1}  \\
\Theta^q_{\lambda^2 \lambda^3}(q^{-\frac{1}{2}} Z_{c,2}^{-1}
Z_y^{-1} Z_{d,1}^{-1})
& = q^{-\frac{1}{2}} Z_{e,3}^{-1} Z_y^{-1} Z_{x,2}^{-1} q^{-\frac{1}{4}} Z_{x,1}^{-1} Z_{a,1}^{-1} \\
& = q^{-\frac{3}{4}} Z_{e,3}^{-1} Z_y^{-1} Z_x^{-1} Z_{a,1}^{-1}\\
\end{split}
\end{equation*}

Thus
$$\Theta^q_{\lambda^{0} \lambda^{1}} \circ \Theta^q_{\lambda^{1} \lambda^{2}} \circ \Theta^q_{\lambda^{2} \lambda^{3}} \circ
\Theta^q_{\lambda^{3} \lambda^{4}} \circ
\Theta^q_{\lambda^{4}
\lambda_{(5)}}(q^{-\frac{1}{4}}Z_{d,1}^{-1} Z_{e,1}^{-1})$$
$$ = Id (q^{-\frac{1}{4}}Z_{d,1}^{-1} Z_{e,1}^{-1})$$.

Now we will consider $\beta$. If we let $\beta$ be as depicted in
Figure~\ref{fig:PentRelSquares} and label the edges of the
pentagons also as represented in Figure~\ref{fig:PentRelSquares} then we
really only need to look at where $Z_{a,3}^{-1} Z_x^{-1}
Z_{c,2}^{-1}$ is mapped to. If we apply the definition of
$\Theta^q_{\lambda \hat{\lambda}}$ and use
Lemma~\ref{Lem:InverseMap} on this monomial we obtain:
\begin{equation*}
\begin{split}
\Theta^q_{\lambda^4 \lambda^0}&(Z_{a,3}^{-1} Z_x^{-1} Z_{c,2}^{-1})  = Z_{c,3}^{-1} (Z_y + Z_y^{-1})^{-1}Z_{e,2}^{-1} \\
\Theta^q_{\lambda^3 \lambda^4}&(Z_{c,3}^{-1} (Z_y +
Z_y^{-1})^{-1}Z_{e,2}^{-1})
 = \Theta^q_{\lambda^3 \lambda^4}(q^{-\frac{1}{4}}Z_{c,3}^{-1}z_{y,3}^{-1} q^{\frac{1}{4}} z_{y,2}^{-1} (1 + Z_y^{-2})^{-1}Z_{e,2}^{-1}) \\
& = q^{-\frac{1}{4}} Z_{e,1}^{-1} Z_y^{-1}Z_x^{-1}(1+ Z_x^{-2}(1+q Z_y^{-2}) )^{-1} Z_{b,3}^{-1}\\
\Theta^q_{\lambda^1 \lambda^0}&(Z_{a,3}^{-1} Z_x^{-1} Z_{c,2}^{-1})  = \Theta^q_{\lambda^1 \lambda^0}(q^{\frac{1}{4}}Z_{a,3}^{-1} Z_{x,3}^{-1} q^{-\frac{1}{4}}Z_{x,2}^{-1}Z_{c,2}^{-1}) \\
& = q^{-\frac{1}{4}} Z_{d,3}^{-1} Z_y^{-1} Z_x^{-1} Z_{a,1}^{-1} \\
 \Theta^q_{\lambda^2 \lambda^1}&(q^{-\frac{1}{4}} Z_{d,3}^{-1} Z_y^{-1} Z_x^{-1} Z_{a,1}^{-1})  =
\Theta^q_{\lambda^2 \lambda^1}(q^{-\frac{1}{4}} Z_{d,3}^{-1} Z_y^{-1} Z_{x,2}^{-1}Z_{x,1}^{-1} Z_{a,1}^{-1}) \\
 &= q^{-\frac{1}{4}} Z_{b,1}^{-1} (Z_x + Z_x^{-1})^{-1} Z_y^{-1}  Z_{d,2}^{-1} \\
\Theta^q_{\lambda^3 \lambda^2}&(q^{-\frac{1}{4}} Z_{b,1}^{-1} (Z_x
+
Z_x^{-1})^{-1} Z_y^{-1}  Z_{d,2}^{-1}) \\
& =
\Theta^q_{\lambda^3 \lambda^2}(q^{-\frac{1}{4}} Z_{b,1}^{-1} (1 + Z_x^{-2})^{-1} Z_{x,1}^{-1} Z_{x,3}^{-1} Z_y^{-1}  Z_{d,2}^{-1}) \\
& = q^{-\frac{1}{4}} Z_{e, 1}^{-1} (1 + Z_y^{-2}(1 + q Z_x^2)^{-1})^{-1} Z_y^{-1} (Z_x + Z_x^{-1})^{-1} Z_{b,3}^{-1} \\
& = q^{-\frac{1}{4}} Z_{e, 1}^{-1} (1 + q^{-1} Z_x^2 Z_y^{-2})^{-1} Z_y^{-1} Z_x  Z_{b,3}^{-1} \\
& = q^{-\frac{1}{4}} Z_{e, 1}^{-1} (Z_x^{-2} + q^{-1} + Z_x^{-2}
Z_y^{-2})^{-1} Z_x^{-2} Z_y^{-1} Z_x  Z_{b,3}^{-1}\\
& = q^{-\frac{1}{4}} Z_{e,1}^{-1} Z_y^{-1}Z_x^{-1}(1+ Z_x^{-2}(1+q
Z_y^{-2}) )^{-1} Z_{b,3}^{-1}
\end{split}
\end{equation*}

Thus
$$\Theta^q_{\lambda^{0} \lambda^{1}} \circ \Theta^q_{\lambda^{1} \lambda^{2}} \circ \Theta^q_{\lambda^{2} \lambda^{3}} \circ
\Theta^q_{\lambda^{3} \lambda^{4}} \circ
\Theta^q_{\lambda^{4} \lambda^{5}}(Z_{a,3}^{-1} Z_x^{-1}
Z_{c,2}^{-1})$$
$$ = Id (Z_{a,3}^{-1} Z_x^{-1} Z_{c,2}^{-1})$$.

\end{proof}

The diagonal exchanges and edge reindexings satisfy the following
relations:

\begin{enumerate}
\item \emph{The Composition Relation}: if $\delta$ and $\gamma$
are each either a diagonal exchange or a edge reindexing then $(\delta
\gamma)(\lambda)=\delta \circ \gamma (\lambda)$.
\item \emph{The Reflexivity Relation}: if $\Delta_i$ is an ith-diagonal
exchange map then $\Delta_i^2(\lambda)=\lambda$.
\item \emph{The Reindexing Relation}: If $\gamma \in S_n$ is a reindexing and $\Delta_i$ is an ith-diagonal exchange map then $\Delta_i \circ \gamma = \gamma \circ \Delta_{\gamma(i)}$.
\item \emph{The Distant Commutativity Relation}: If $\lambda_i$ and $\lambda_j$ are edges of the ideal triangulation $\lambda \in \Lambda(S)$ that do not belong to the same triangle then $\Delta_i \circ \Delta_j(\lambda) = \Delta_j \circ \Delta_i(\lambda)$.
\end{enumerate}

We now state the following two results of Penner. Refer to  ~\cite{Pen1} for
their proofs.

\begin{theorem}
\label{thm:Penner1}
Given two ideal triangulations $\lambda$, $\hat{\lambda}$, there
exists a finite sequence of ideal triangulations
$\lambda=\lambda^0$, $\lambda^1$, ...,$\lambda^m=\hat{\lambda}$
such that $\lambda^{k+1}$ is obtained from $\lambda^k$ by a single
diagonal exchange or by edge reindexing.
\end{theorem}
\qed

\begin{theorem}
\label{thm:Penner2}
Given two ideal triangulations $\lambda$, $\hat{\lambda}$ and two
sequences of ideal triangulations $\lambda=\lambda^0$,
$\lambda^1$, ...,$\lambda^m=\hat{\lambda}$ and
$\lambda=\hat{\lambda}^0$, $\hat{\lambda}^1$,
...,$\hat{\lambda}^m=\hat{\lambda}$ such that $\lambda^{k+1}$ is
obtained from $\lambda^k$ by a single diagonal exchange or by edge
reindexing and $\hat{\lambda}^{k+1}$ is obtained from
$\hat{\lambda}^k$ by a single diagonal exchange or by edge
reindexing, these two sequences can be related to each other by
applications of the following moves and their inverses:

\begin{enumerate}
\item Use the \emph{The Composition Relation} to replace $...,\lambda^{k}, \delta(\lambda^{k}), \gamma \circ
    \delta(\lambda^{k}),...$ with $...,\lambda^{k}, (\gamma \delta)(\lambda^{k}),...$ where $\delta$ and $\gamma$
are each either a diagonal exchange or a edge reindexing.
\item Use \emph{The Reflexivity Relation} to replace $...,\lambda^{k},...$ with
    $...,\lambda^{k},\Delta_i(\lambda^{k}),\lambda^{k}...$
\item Use \emph{The Reindexing Relation} to replace $...,\lambda^{k}, \gamma(\lambda^{k}),
    \Delta_i(\gamma(\lambda^{k})),...$ with $...,\lambda^{k}, \Delta_{\gamma(i)}(\lambda^{k}),
    \gamma(\Delta_{\gamma(i)}(\lambda^{k})),...$ where $\gamma \in S_n$ is an edge reindexing.
\item Use \emph{The Distant Commutativity Relation} to replace $...,\lambda^{k},...$ with
    $...,\lambda^{k},\Delta_i(\lambda^{k}),\Delta_j(\Delta_i(\lambda^{k})),\Delta_j(\lambda^{k}),\lambda^{k}...$ where
    $\lambda_i$ and $\lambda_j$ are two edges of $\lambda^{k}$ that do not belong to the same triangle.
\item Use \emph{The Pentagon Relation} to replace $...,\lambda^{k},...$ with $...,\lambda^{k},\Delta_i(\lambda^{k}),
    \Delta_j \circ \Delta_i(\lambda^{k}), \Delta_i \circ \Delta_j \circ \Delta_i(\lambda^{k}),\Delta_j \circ \Delta_i
    \circ \Delta_j \circ \Delta_i(\lambda^{k}),\alpha_{i \leftrightarrow j}(\lambda^{k}), \lambda^{k},...$ where
    $\lambda_i$ and $\lambda_j$ are two diagonals of a pentagon of $\lambda^{k}$.
\end{enumerate}

\end{theorem}
\qed

Note: If we are given two ideal triangulations $\lambda$ and $\hat{\lambda}$ we can find a sequence of ideal
triangulations $\lambda = \lambda^{0}$, $\lambda^{1}$, ...,$\lambda^{m}=\hat{\lambda}$
where each $\lambda^{k+1}$ is obtained from $\lambda^{k}$ by a
diagonal exchange or by an edge reindexing. Define  $\Theta^q_{\lambda
\hat{\lambda}}$, as the composition of the
$\Theta^q_{\lambda^{k} \lambda^{k+1}}$. Lemma~\ref{Lem:InverseMap} and Lemma~\ref{lem:PentagonRelation}
 along with Theorem ~\ref{thm:Penner1} and Theorem ~\ref{thm:Penner2} show that this $\Theta^q_{\lambda^{k} \lambda^{k+1}}$ is
independent of the choice of the sequence of $\lambda^{k}$.

\begin{theorem}
Given ideal triangulations $\lambda$, $\lambda'$, $\lambda''$ then,
\begin{equation*}
\Theta^q_{\lambda \lambda''} = \Theta^q_{\lambda \lambda'} \circ
\Theta^q_{\lambda' \lambda''}
\end{equation*}
\end{theorem}
\begin{proof}

This result simply follows from the definition.

\end{proof}

\section{Punctured Tori}

\begin{figure}[htb]
\centering
 \SetLabels
\endSetLabels
\centerline{\AffixLabels{
\includegraphics{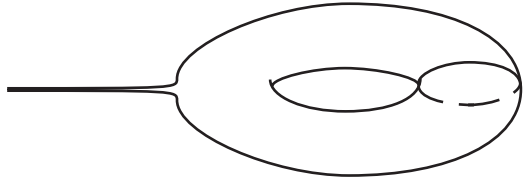}}}
\caption{The Once Punctured Torus}
\end{figure}

\begin{figure}[htb]
\centering
\SetLabels
\endSetLabels
\centerline{\AffixLabels{
\includegraphics{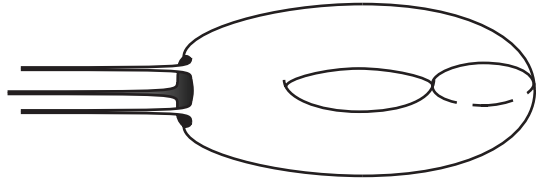}}}
\caption{A Torus With a Wide
Hole and $p \geq 1$ Punctures}
\end{figure}

Let the \emph{once-punctured torus} be the surface obtained by
removing one point from a torus. Let a \emph{torus with a wide
hole and $p \geq 1$ punctures} be the surface that is obtained
from the compact surface of genus one with one boundary component
by removing $p \geq 1$ punctures from its boundary but none from
its interior.

In this section, $S$ will denote either a once punctured torus or
a torus with a wide hole and $p \geq 1$ punctures. Let $\Sigma(S)$
be the set of simple closed unoriented curves in $S$

\begin{definition}
For an ideal triangulation $\lambda$ of $S$ with edges
$\lambda_1$,...,$\lambda_n$ and given $\alpha \in \Sigma(S)$,
define $\alpha$ as \emph{$\lambda$-simple} if it meets each
$\lambda_i$ in at most one point.
\end{definition}

\begin{definition}
Given an ideal triangulation $\lambda$ of $S$ with edges
$\lambda_1$,...,$\lambda_n$ and given $\alpha \in \Sigma(S)$,
define $\lambda$ as \emph{$\alpha$-simple} if and only if $\alpha$
is $\lambda$-simple.
\end{definition}

\begin{figure}[htb]
\centering
 \SetLabels
( 0 * .24) $\lambda_j$\\
( 0.85 * .85) $\lambda_j$\\
  ( .2 * .91) $\lambda_k$\\
  ( 0.5 * -0.05) $\lambda_k$\\
  ( .26 * .1) $\alpha$\\
  ( .52 * .9) $\alpha$\\
   ( .42 * .22) $T_1$\\
  ( .48 * .7) $T_2$\\
\endSetLabels
\centerline{\AffixLabels{
\includegraphics{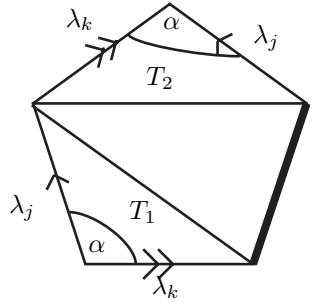}}}
\caption{Two Dimensional Torus With a Wide Hole}
\label{fig:QuadWithBoundaryTorus}
\end{figure}

We may now state the two main theorems of this section.

\begin{theorem}
\label{thm:TorusWideHole}
 Let $S$ be either a once punctured torus or a torus with a wide hole and $p \geq 1$ punctures.
 There exists a family of $T_\alpha^\lambda \in \mathcal{T}_\lambda^{q^{\frac{1}{4}}}(\alpha)$, with $\lambda$
ranging over all ideal triangulations of $S$ and $\alpha$ over all
essential simple closed curves of $S$, which satisfies:

\begin{enumerate}
    \item If $\alpha$ is in $\Sigma(S)$ and $\lambda$ and
$\hat{\lambda}$ are two triangulations of $S$, then $
\Theta_{\lambda
\hat{\lambda}}(T_\alpha^{\hat{\lambda}})=T_\alpha^\lambda$.

    \item as $q \rightarrow 1$, $T_{\alpha}^{\lambda}$ converges to the non-quantum
    trace function $T_{\alpha}$ in $\mathcal{T}(S)$

    \item If $\alpha$ and $\beta$ are disjoint, $T_{\alpha}^{\lambda}$ and
    $T_{\beta}^{\lambda}$ commute.
    \item If $\alpha$ meets each edge of $\lambda$ at most once then $T_{\alpha}^{\lambda}$
    is obtained from the classical trace $T_{\alpha}$ of Section 2.3 by multiplying each
    monomial by the Weyl ordering coefficient.
\end{enumerate}
\end{theorem}

\begin{figure}[htb]
\centering
 \SetLabels
 ( 0 * .2) $\alpha$\\
  ( .2 * .2) $\beta$\\
  ( .5 * .2) $\alpha\beta$\\
  ( .9 * .06) $\beta\alpha$\\
\endSetLabels
\centerline{\AffixLabels{
\includegraphics{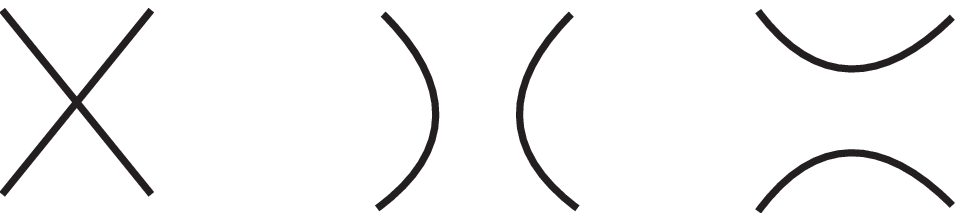}}}
\caption{Resolving Crossing in the Torus} \label{fig:SkeinTorus}
\end{figure}

\begin{theorem}
\label{thm:TorusWideHole2} The traces $T_{\alpha}^{\lambda}$ of
Theorem~\ref{thm:TorusWideHole} satisfy the following property: If
$\alpha$ and $\beta$ meet in one point, and if
    $\alpha \beta$ and $\beta \alpha$ are obtained by resolving
    the intersection point as in Figure~\ref{fig:SkeinTorus}, then
    $$T_{\alpha}^{\lambda}
    T_{\beta}^{\lambda} = q T_{\alpha \beta}^{\lambda} + q^{-1} T_{\beta \alpha}^{\lambda}.$$

    In addition,
$T_{\alpha}^{\lambda}$ with $\alpha$ non-separating is the only
one which satisfy this property and conditions (1) and (4) of
Theorem~\ref{thm:TorusWideHole}.
\end{theorem}

We restrict our attention to the case where $S$ is a torus with
a wide hole and $p \geq 1$ punctures. The case of the
once-punctured torus is similar, and simpler.

A key result used to prove Theorem~\ref{thm:TorusWideHole} is the
following proposition.

\begin{proposition}
\label{prop:KeyPropertyTorus}
 For an essential curve $\alpha$ in $\Sigma(S)$ and $\alpha$-simple ideal triangulations $\lambda$ and $\hat{\lambda}$ of $S$ there exists a sequence of ideal
 triangulations $\lambda=\lambda^0$ , $\lambda^1$ ,$\lambda^2$
 ,..., $\lambda^{m-1}$
,$\lambda^m =\hat{\lambda}$ such that each $\lambda^{i+1}$ is
obtained from $\lambda^i$ by a diagonal exchange and $\lambda^i$
is $\alpha$-simple.

\end{proposition}

\begin{proof}
To prove Proposition~\ref{prop:KeyPropertyTorus} we must prove the
following lemmas.

\begin{lemma}
\label{lem:KeyPropertyTorusLemma}
    If $\lambda$ and $\hat{\lambda}$ are $\alpha$-simple ideal triangulations of $S$, there exists a sequence of
    $\alpha$-simple ideal triangulations $\lambda=\lambda^0$ ,
    $\lambda^1$ ,..., $\lambda^m$ and a sequence of $\alpha$-simple ideal triangulations $\hat{\lambda}=\hat{\lambda}^0$ ,
    $\hat{\lambda}^1$ ,..., $\hat{\lambda}^n$ such that:

\begin{enumerate}
    \item $\lambda^{l+1}$ is obtained from $\lambda^l$ by a diagonal exchange
    and $\hat{\lambda}^{l+1}$ is obtained from $\hat{\lambda}^l$ by a diagonal
    exchange
     \item there exists edges $\lambda_i$ and $\hat{\lambda}_j$ of
 $\lambda^m$ and $\hat{\lambda}^n$, respectively, such that one component $C_1$ of $S-\lambda_i$ and one component $\hat{C}_1$ of
    $S-\hat{\lambda}_j$ are tori with exactly one spike at infinity.
    \item $\lambda^m$ and $\hat{\lambda}^n$ coincide outside of $C_1$.
 \end{enumerate}
\end{lemma}

\begin{figure}[htb]
\centering
 \SetLabels
 ( .5 * 1.04) $S_1$\\
  ( 1.06 * .5) $S_2$\\
  ( .5 * -0.08) $S_3$\\
  ( -0.06 * .5) $S_4$\\
  ( -0.06 * .85) $S_5$\\
\endSetLabels
\centerline{\AffixLabels{
\includegraphics{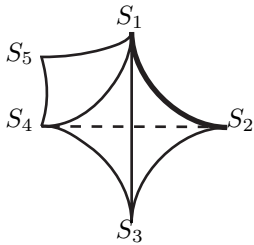}}}
\caption{Quadrilateral in the Surface} \label{fig:QuadTorus}
\end{figure}

\begin{proof}
    We will use the following proof for both $\lambda$ and $\hat{\lambda}$.
    Consider a quadrilateral as represented in Figure~\ref{fig:QuadTorus}
    with vertices at infinity $S_1,S_2,S_3,S_4$, occurring counterclockwise in this order, and where  and the edge from $S_1$ to
    $S_2$ is a boundary curve, the edges in this quadrilateral
    are edges $S_1$ to $S_2$, $S_1$ to $S_3$, $S_1$ to $S_4$, $S_2$ to $S_3$ and $S_3$ to
    $S_4$. If $S_1=S_2$ then we are done. Assume from now on that
    $S_1 \neq S_2$

    \textbf{Case 1:}  If $S_3=S_1$ and $S_4=S_1$ then doing a diagonal
    exchange on both $\lambda$ and $\hat{\lambda}$ in this quadrilateral lowers the number of edges
    ending at $S_1$ by one. Also since $\alpha$ cannot cross the edge connecting $S_1$ to $S_2$
    when you do this diagonal exchange the resulting triangulation remains $\alpha$-simple.

    \textbf{Case 2:}  If $S_1=S_3$ and $S_1 \neq S_4$, then doing a diagonal
    exchange on both $\lambda$ and $\hat{\lambda}$ in this quadrilateral lowers the number of edges
    ending at $S_1$ by two. Also, the resulting triangulation remains $\alpha$-simple.

    \textbf{Case 3:} If  $S_1 \neq S_3$ and $S_1 \neq S_4$, then doing a diagonal exchange on both $\lambda$ and $\hat{\lambda}$ in
    this quadrilateral decreases the number of edges ending at
    $S_1$ by one. Also, the resulting triangulation remains $\alpha$-simple.

    \textbf{Case 4:} If  $S_1 \neq S_3$ and $S_1=S_4$. then after doing a diagonal exchange on both $\lambda$ and $\hat{\lambda}$ in this
    quadrilateral if we consider the new quadrilateral created
    with edges $S_1$ to $S_2$, $S_2$ to $S_4$, $S_1$ to $S_4$ and
    a new point $S_5$ and edges $S_1$ to $S_5$ and $S_4$ to $S_5$
    then we see that we are again in Case 2 or Case 3. Thus after another
    diagonal exchange on both $\lambda$ and $\hat{\lambda}$ in this new quadrilateral we reduce the
    number of edges ending at $S_1$ by one or two. For the same reason as above after the first diagonal exchange the
    resulting triangulation remains $\alpha$-simple and similarly after the second diagonal exchange the resulting
    triangulation remains $\alpha$-simple.

    Now if we repeat this process until there are only two edges going to $S_1$ then we can
    effectively ``forget'' about the point at infinity $S_1$ and then repeat this process for another
    point at infinity. By the method of the proof we automatically get $\lambda^m$ and $\hat{\lambda}^n$ to coincide outside of $C_1$.
\end{proof}

\begin{figure}[htb]
\centering
 \SetLabels
 (0.9 * 0.16) $T_1$\\
 (0.17 * 0.76) $T_1$\\
 (0.835 * 0.07) $\alpha$\\
 (0.47 * 0.22) $\alpha$\\
 (0.09 * 0.28) $\alpha$\\
 (0.95 * 0.83) $\alpha$\\
 (0.51 * 0.88) $\alpha$\\
 (0.09 * 0.68) $\alpha$\\
 (0.26 * 0.7) $\lambda_i$\\
( .32 * 0.8) $=$\\
  ( 1.06 * .8) $=$\\
  ( .68 * 0.2) $=$\\
  ( 0.275 * .28) \tiny{$1$}\\
  ( 0.34 * .28) \tiny{$Diagonal$}\\
  ( 0.33 * .16) \tiny{$Exchange$}\\
  ( 0.64 * .86) \tiny{$2$}\\
  ( 0.7* .86) \tiny{$Diagonal$}\\
  ( 0.69 * .76) \tiny{$Exchanges$}\\
\endSetLabels
\centerline{\AffixLabels{
\includegraphics{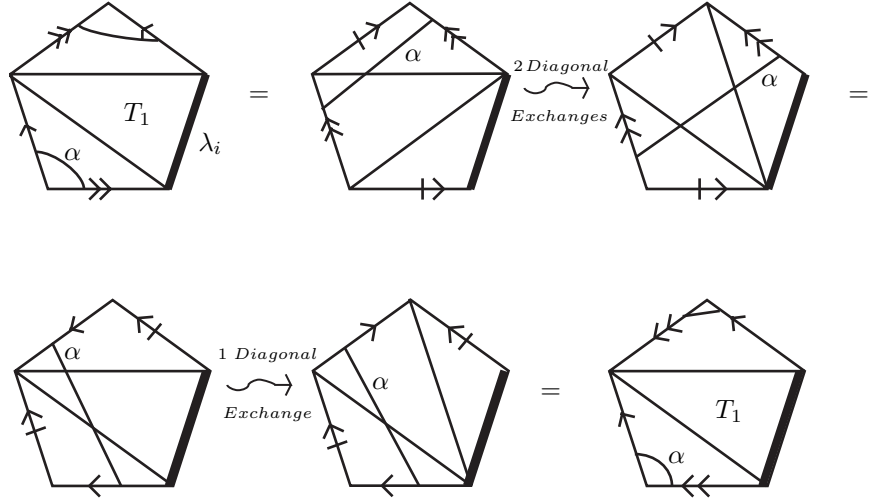}}}
\caption{Pentagon Moves} \label{fig:QuadTorus2D}
\end{figure}

\begin{figure}[htb]
\centering
 \SetLabels
\endSetLabels
\centerline{\AffixLabels{
\includegraphics{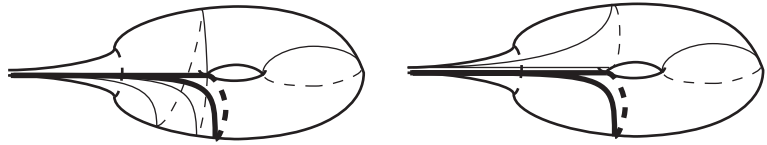}}}
\caption{Torus Slope Changes} \label{fig:QuadTorus3D}
\end{figure}

\begin{lemma}
\label{lem:TriangleTorus} After changing the ideal triangulation
$\lambda$ by $\alpha$-simple diagonal exchanges, we can
arrange that, for the component $C_1$ of $S-\lambda_i$ which is a
torus with one spike at infinity, the triangle $T_1$ containing
$\partial C_1$ is disjoint from $\alpha$.

\end{lemma}

\begin{proof}
This can be accomplished using moves similar to the ones
represented in Figure~\ref{fig:QuadTorus2D}.

\end{proof}

Change both triangulations $\lambda$ and $\hat{\lambda}$ such that
the triangle $T_1$ containing $\partial C_1$ and
    $\partial \hat{C}_1$ is disjoint from $\alpha$ as in Lemma~\ref{lem:TriangleTorus}. If we consider the
    two nonboundary
    edges $\lambda_1$, $\lambda_2$ that make up the triangle $T_1$, they are completely determined
    by how many times they wrap around the boundary. Also, if edge $\lambda_1$ wraps around the boundary $k$ times then
    edge $\lambda_2$ is restricted to wrap around the boundary $k+1$ or $k-1$ times.  Now it is clear from Figure~\ref{fig:QuadTorus2D}
    that if edges $\lambda_1$ and $\lambda_2$ wrap around the boundary $k$ and $k+1$ times respectively that through a series of
    $\alpha$-simple diagonal exchanges we can move to ideal triangulation $\lambda'$ where the two edges $\lambda_1'$ and $\lambda_2'$,
    which are edges in the triangle on the boundary, wrap around the boundary either $k+1$ and $k+2$ time respectively or $k-1$ and $k$ times
    respectively. This is further illustrated in Figure~\ref{fig:QuadTorus3D}. Thus we can change $\lambda$ and $\hat{\lambda}$ so that the
    edges $\lambda_1$ and $\lambda_2$ which make up the triangle on the boundary of $C_1$ and $\hat{C}_1$ coincide.

\begin{figure}[htb]
\centering
 \SetLabels
\endSetLabels
\centerline{\AffixLabels{
\includegraphics{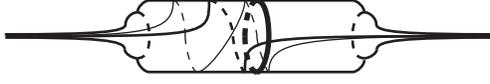}}}
\caption{Cylinder in the Torus} \label{fig:CylinderTorus3D}
\end{figure}

If we cut the surface along the edges $\lambda_1$ and $\lambda_2$
from above we are left with a cylinder with two spikes at
infinity, four edges going between those spikes and $\alpha$ as a
meridian of the cylinder as depicted in Figure~\ref{fig:CylinderTorus3D}.
It is clear that, if we perform a diagonal exchange along any one
of the two edges that spiral around the cylinder, the resulting
ideal triangulation remains $\alpha$-simple. In addition, the
diagonal exchanges of this type enable us to go between any two
ideal triangulations of the cylinder.

Thus, we can always reduce an
ideal triangulation to the case
    with only one spike at infinity.

    This concludes the proof of Proposition ~\ref{prop:KeyPropertyTorus}
\end{proof}

We now prove Theorem~\ref{thm:TorusWideHole} which we restate for
the sake of the reader.

\begin{theorem}
\label{thm:TorusWideHole1}
 There exists a family of $T_\alpha^\lambda \in \mathcal{T}_\lambda^{q^{\frac{1}{4}}}(\alpha)$ with $\lambda$
ranging over all ideal triangulations of $S$ and $\alpha$ over all essential
simple closed curves of $S$, which satisfies:

\begin{enumerate}
    \item If $\alpha$ is in $\Sigma(S)$ and $\lambda$ and
$\hat{\lambda}$ are two triangulations of $S$, then $
\Theta_{\lambda
\hat{\lambda}}(T_\alpha^{\hat{\lambda}})=T_\alpha^\lambda.$

    \item as $q \rightarrow 1$, $T_{\alpha}^{\lambda}$ converges to the non-quantum
    trace function $T_{\alpha}$ in $\mathcal{T}(S).$

    \item If $\alpha$ and $\beta$ are disjoint, $T_{\alpha}^{\lambda}$ and
    $T_{\beta}^{\lambda}$ commute.
    \item If $\alpha$ meets each edge of $\lambda$ at most once then $T_{\alpha}^{\lambda}$
    is obtained from the classical trace $T_{\alpha}$ of Section 2.3 by multiplying each
    monomial by the Weyl ordering coefficient.
\end{enumerate}
\end{theorem}

\begin{proof}

For an ideal triangulation $\lambda$ and a simple closed curve
$\alpha$ in $S$ which is $\lambda$-simple, define
$T_{\alpha}^{\lambda}$ to be obtained from the non-quantum trace
function $T_{\alpha}$ by multiplying each monomial of $T_{\alpha}$
with the Weyl quantum ordering coefficient.

When $\alpha$ is not $\lambda$-simple and is not homotopic to the boundary, one easily finds an
$\alpha$-simple ideal triangulation $\lambda^*$. In this case define
$T_{\alpha}^{\lambda}=\Theta_{\lambda
\lambda^*}(T_{\alpha}^{\lambda^*})$ where $T_{\alpha}^{\lambda^*}
\in \hat{\mathcal{T}}_{\lambda^*}^{\frac{1}{4}}(\alpha)$ is
defined by the previous case.

When $\alpha$ is not $\lambda$-simple and is homotopic to the boundary define
$T_{\alpha}^{\lambda}$ to be obtained from the non-quantum trace
function $T_{\alpha}$ by multiplying each monomial of $T_{\alpha}$
with the Weyl quantum ordering coefficient.

In the case where $\alpha$ is not homotopic to the boundary let us show that $T_{\alpha}^{\lambda}$ is well defined,
namely is independent of the choice of the $\alpha$-simple ideal
triangulation $\lambda^*$. The main step is to prove the following lemma.

\begin{lemma}
\label{lem:TorusWideHoleIndepCoord}
 For any two $\alpha$-simple ideal triangulations $\lambda , \hat{\lambda}$ of $S$,  $
\Theta_{\lambda
\hat{\lambda}}(T_\alpha^{\hat{\lambda}})=T_\alpha^\lambda$.

\end{lemma}
\begin{proof}
The main step in the proof is the following.

\begin{figure}[htb]
\centering
 \SetLabels
 ( 0.2 * 1.06) $\lambda_j$\\
  ( .8 * 1.06) $\hat{\lambda}_j$\\
  ( .2 * 0.06) $\lambda_l$\\
  ( .8 * .08) $\hat{\lambda}_l$\\
( .16 * 0.56) $\lambda_i$\\
  ( .84 * .56) $\hat{\lambda}_i$\\
  ( -0.04 * 0.5) $\lambda_m$\\
  ( 1.04 * 0.5) $\hat{\lambda}_k$\\
  ( .4 * 0.5) $\lambda_k$\\
  ( .62 * .5) $\hat{\lambda}_m$\\
   ( .11 * 0.75) $T_1$\\
  ( .25 * .3) $T_2$\\
   ( .89 * 0.75) $T_1$\\
  ( .75 * .3) $T_2$\\
\endSetLabels
\centerline{\AffixLabels{
\includegraphics{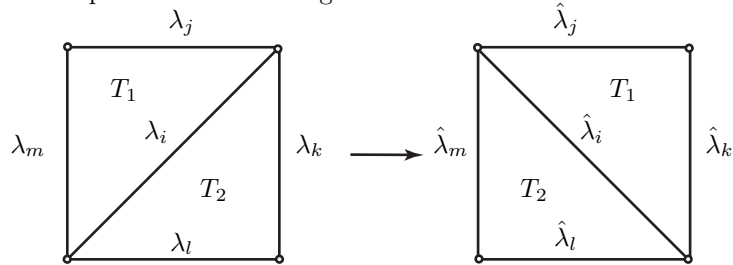}}}
\caption{Diagonal Exchange} \label{Fig:DiagExchangeTorus}
\end{figure}


\begin{lemma}
\label{lem:TorusWideHoleFirst}
 Given two $\alpha$-simple ideal triangulations $\lambda$ and $\hat{\lambda}$ which differ by only one diagonal
exchange,  $ \Theta_{\lambda
\hat{\lambda}}(T_\alpha^{\hat{\lambda}})=T_\alpha^\lambda$.
\end{lemma}

\begin{proof}
We must separate this into four cases, each represented in Figure~\ref{Fig:DiagExchangeTorusCase1}. In each case we have
$\alpha \in \Sigma(S)$ and $\alpha$-simple ideal triangulations $\lambda$ and
$\hat{\lambda}$ with $\lambda$ and $\hat{\lambda}$ differing
only by exchanging edge $\lambda_i$ and $\hat{\lambda}_i$. The
edges of the quadrilaterals involved in the diagonal exchange are
$\lambda_i$, $\lambda_j$, $\lambda_k$, $\lambda_l$, $\lambda_m$,
$\hat{\lambda}_i$, $\hat{\lambda}_j$, $\hat{\lambda}_k$,
$\hat{\lambda}_l$ and $\hat{\lambda}_m$ with the quadrilaterals
edges in $\lambda$ and $\hat{\lambda}$ ordered as represented in Figure~\ref{Fig:DiagExchangeTorus}. This
lemma follows from simple calculations of which we will do one of
the four possible cases.

\vspace{2mm}

\begin{figure}[htb]
\centering
 \SetLabels
 ( 0.08 * 1.02) $\lambda_j$\\
  ( .34 * 1.02) $\hat{\lambda}_j$\\
  ( 0.08 * 0.63) $\lambda_l$\\
  ( .34 * .62) $\hat{\lambda}_l$\\
( .04 * 0.795) $\lambda_i$\\
  ( .33 * .795) $\hat{\lambda}_i$\\
  ( -0.02 * 0.83) $\lambda_m$\\
  ( 0.255 * 0.83) $\hat{\lambda}_m$\\
  ( .17 * 0.83) $\lambda_k$\\
  ( .44 * .83) $\hat{\lambda}_k$\\
  ( .11 * 0.78) $\alpha$\\
  ( .39 * .86) $\alpha$\\
  (0.21 * 0.54) Case 1\\
  ( 0.66 * 1.02) $\lambda_j$\\
  ( .92 * 1.02) $\hat{\lambda}_j$\\
  ( 0.66 * 0.63) $\lambda_l$\\
  ( .92 * .62) $\hat{\lambda}_l$\\
( .62 * 0.795) $\lambda_i$\\
  ( .91 * .795) $\hat{\lambda}_i$\\
  ( 0.56 * 0.83) $\lambda_m$\\
  ( 0.835 * 0.83) $\hat{\lambda}_m$\\
  ( .755 * 0.83) $\lambda_k$\\
  ( 1.02 * .83) $\hat{\lambda}_k$\\
  ( .63 * 0.89) $\alpha$\\
  ( .95 * .89) $\alpha$\\
  (0.79 * 0.54) Case 2\\
  ( 0.08 * .33) $\lambda_j$\\
  ( .34 * 0.33) $\hat{\lambda}_j$\\
  ( 0.08 * -0.06) $\lambda_l$\\
  ( .34 * -0.06) $\hat{\lambda}_l$\\
( .04 * 0.105) $\lambda_i$\\
  ( .33 * .105) $\hat{\lambda}_i$\\
  ( -0.02 * 0.14) $\lambda_m$\\
  ( 0.255 * 0.14) $\hat{\lambda}_m$\\
  ( .17 * 0.14) $\lambda_k$\\
  ( .44 * .14) $\hat{\lambda}_k$\\
  ( .06 * 0.24) $\alpha$\\
  ( .33 * .24) $\alpha$\\
  (0.21 * -0.1) Case 3\\
  ( 0.66 * 0.33) $\lambda_j$\\
  ( 0.92 * 0.33) $\hat{\lambda}_j$\\
  ( 0.66 * -0.06) $\lambda_l$\\
  ( .92 * -.07) $\hat{\lambda}_l$\\
( .62 * 0.105) $\lambda_i$\\
  ( .91 * .105) $\hat{\lambda}_i$\\
  ( 0.56 * 0.14) $\lambda_m$\\
  ( 0.835 * 0.14) $\hat{\lambda}_m$\\
  ( .755 * 0.14) $\lambda_k$\\
  ( 1.02 * .14) $\hat{\lambda}_k$\\
  ( .68 * 0.24) $\alpha$\\
  ( .94 * .24) $\alpha$\\
  (0.79 * -0.1) Case 4\\
\endSetLabels
\centerline{\AffixLabels{
\includegraphics{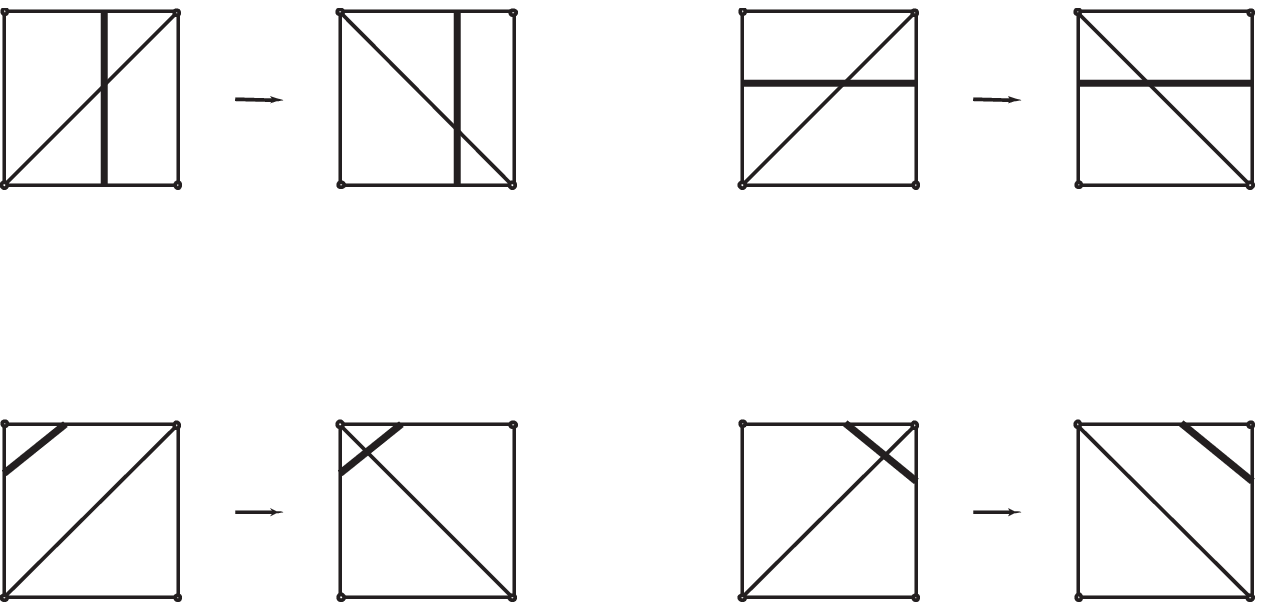}}}
\caption{} \label{Fig:DiagExchangeTorusCase1}
\end{figure}

\textbf{Case 1:} Let all edges of the quadrilateral be distinct
and let the diagonal $\lambda_i$ go from the vertex adjoining
$\hat{\lambda}_j$ to $\hat{\lambda}_k$ to the vertex adjoining
$\hat{\lambda}_l$ to $\hat{\lambda}_m$. Also, let $\alpha$ cross
edges $\hat{\lambda}_j$, $\hat{\lambda}_i$, and $\hat{\lambda}_l$
as depicted in Case 1 of Figure~\ref{Fig:DiagExchangeTorusCase1}. Then by simple
calculations we have:

\begin{equation*}
\begin{split}
 T_{\alpha}^{\hat{\lambda}}  = & \hat{A}(\hat{Z}_{j,2} \hat{Z}_i \hat{Z}_{l,1} + \hat{Z}_{j,2}^{-1} \hat{Z}_i \hat{Z}_{l,1}^{-1} +
 \hat{Z}_{j,2}^{-1} \hat{Z}_i^{-1} \hat{Z}_{l,1}^{-1})\hat{B} +  \hat{C}(q^{-\frac{1}{2}}\hat{Z}_{j,2} \hat{Z}_i \hat{Z}_{l,1}^{-1})\hat{D} +
 \hat{E}(q^{\frac{1}{2}}\hat{Z}_{j,2}^{-1} \hat{Z}_i \hat{Z}_{l,1})\hat{F} \\
 & \hat{A},\hat{B},\hat{C},\hat{D},\hat{E},\hat{F}  \in \mathcal{T}_{\hat{\lambda}}^{q^{\frac{1}{4}}}\\
\end{split}
\end{equation*}
and
 \begin{equation*}
\begin{split}
 \Theta_{\lambda \hat{\lambda}}(T_{\alpha}^{\hat{\lambda}})&=\Theta_{\lambda \hat{\lambda}}(\hat{A}(\hat{Z}_{j,2} \hat{Z}_i \hat{Z}_{l,1} + \hat{Z}_{j,2}^{-1} \hat{Z}_i \hat{Z}_{l,1}^{-1} +
 \hat{Z}_{j,2}^{-1} \hat{Z}_i^{-1} \hat{Z}_{l,1}^{-1})\hat{B} +  \hat{C}(q^{-\frac{1}{2}}\hat{Z}_{j,2} \hat{Z}_i \hat{Z}_{l,1}^{-1})\hat{D} +
 \hat{E}(q^{-\frac{1}{2}}\hat{Z}_{j,2}^{-1} \hat{Z}_i \hat{Z}_{l,1})\hat{F})\\
 &= A(\Theta_{\lambda \hat{\lambda}}(\hat{Z}_{j,2}^{-1} \hat{Z}_i^{-1} \hat{Z}_{l,1}^{-1}\hat{Z}_{j,2}^2 \hat{Z}_i^2 \hat{Z}_{l,1}^2 +
 q^{-1}\hat{Z}_{j,2}^{-1} \hat{Z}_i^{-1} \hat{Z}_{l,1}^{-1}\hat{Z}_i^2 +
 \hat{Z}_{j,2}^{-1} \hat{Z}_i^{-1} \hat{Z}_{l,1}^{-1}))B\\
 &  \quad \quad \quad + C(\Theta_{\lambda \hat{\lambda}}(q^{\frac{1}{2}}\hat{Z}_{j,2}^{-1} \hat{Z}_i^{-1} \hat{Z}_{l,1}^{-1}\hat{Z}_{j,2}^{2} \hat{Z}_i^{2} ))D
 + E(\Theta_{\lambda \hat{\lambda}}(q^{-\frac{3}{2}}\hat{Z}_{j,2}^{-1} \hat{Z}_i^{-1} \hat{Z}_{l,1}^{-1} \hat{Z}_i^2 \hat{Z}_{l,1}^2))F\\
   & =  A(Z_{j,1}^{-1}(Z_i + Z_i^{-1})^{-1}Z_{l,2}^{-1}(1+ qZ_i^2)Z_{j,1}^2Z_i^{-2}(1+qZ_i^2)Z_{l,2}^2 \\
&\quad \quad \quad + q^{-1} Z_{j,1}^{-1}(Z_i + Z_i^{-1})^{-1} Z_{l,2}^{-1} Z_i^{-2} + Z_{j,1}^{-1}(Z_i + Z_i^{-1})^{-1}Z_{l,2}^{-1})B\\
 & \quad \quad \quad+ C(q^{\frac{1}{2}} Z_{j,1}^{-1} (Z_i + Z_i^{-1})^{-1} Z_{l,2}^{-1} (1 + qZ_i^2) Z_{j,1}^2 Z_i^{-2})D\\
 & \quad \quad \quad+ E(q^{-\frac{3}{2}} Z_{j,1}^{-1} (Z_i + Z_i^{-1})^{-1} Z_{l,2}^{-1} Z_i^{-2} (1+qZ_i^2)Z_{l,2}^2)F\\
& = A(Z_{j,1}^{-1}Z_i(1 + Z_i^2)^{-1}(1+ Z_i^2)Z_{l,2}^{-1}Z_{j,1}^2Z_i^{-2}(1+qZ_i^2)Z_{l,2}^2 \\
& \quad \quad \quad + Z_{j,1}^{-1}(Z_i + Z_i^{-1})^{-1} Z_i^{-2}Z_{l,2}^{-1}  + Z_{j,1}^{-1}(Z_i + Z_i^{-1})^{-1}Z_{l,2}^{-1})B\\
& \quad \quad \quad + C(q^{\frac{1}{2}} Z_{j,1} Z_i^{-1} Z_{l,2}^{-1})D + E(q^{-\frac{1}{2}} Z_{j,1}^{-1} Z_i^{-1} Z_{l,2} )F\\
& = A(Z_{j,1}(Z_i^{-1}+Z_i)Z_{l,2}+ Z_{j,1}^{-1}Z_i^{-1}(1 + Z_i^{-2})^{-1}( Z_i^{-2}  + 1)Z_{l,2}^{-1})B\\
& \quad \quad \quad + C(q^{\frac{1}{2}} Z_{j,1} Z_i^{-1} Z_{l,2}^{-1})D + E(q^{-\frac{1}{2}} Z_{j,1}^{-1} Z_i^{-1} Z_{l,2} )F\\
& = A(Z_{j,1}Z_iZ_{l,2}+Z_{j,1}Z_i^{-1}Z_{l,2}+
Z_{j,1}^{-1}Z_i^{-1}Z_{l,2}^{-1})B\\
  & \quad \quad \quad+ C(q^{\frac{1}{2}} Z_{j,1} Z_i^{-1} Z_{l,2}^{-1})D + E(q^{-\frac{1}{2}} Z_{j,1}^{-1} Z_i^{-1} Z_{l,2} )F = T_{\alpha}^{\lambda}\\
  & \quad \quad \quad \quad \quad \quad   A,B,C,D,E,F  \in \mathcal{T}_{\lambda}^{q^{\frac{1}{4}}}\\
\end{split}
\end{equation*}

\end{proof}

 Lemma~\ref{lem:TorusWideHoleIndepCoord} is now a direct corollary of Proposition~\ref{prop:KeyPropertyTorus} and Lemma~\ref{lem:TorusWideHoleFirst}.
\end{proof}

 Lemma~\ref{lem:TorusWideHoleIndepCoord} proves that the definition of $T_{\alpha}^{\lambda}$ is independent of the choice of
 triangulation $\lambda^*$.


Consider two triangulations $\lambda$ and $\hat{\lambda}$ and an essential simple closed curve $\alpha$.
 From the above
definition we have that there exists an $\alpha$-simple ideal triangulation $\lambda^{*}$
 such that $T_{\alpha}^{\lambda} =
\Theta _{\lambda \lambda^{*}} (T_{\alpha}^{\lambda^{*}} )$ and
$T_{\alpha}^{\hat{\lambda}} = \Theta _{\hat{\lambda} \lambda^{*}}
(T_{\alpha}^{\lambda^{*}} )$. Thus $\Theta_{\lambda
\lambda^{*}}^{-1}(T_{\alpha}^{\lambda}) = \Theta_{\hat{\lambda}
\lambda^{*}}^{-1}(T_{\alpha}^{\hat{\lambda}})$, which implies that
$T_{\alpha}^{\lambda} = \Theta_{\lambda \lambda^{*}} (
\Theta_{\hat{\lambda}
\lambda^{*}}^{-1}(T_{\alpha}^{\hat{\lambda'}}))$. Since we know
that $\Theta_{\hat{\lambda} \lambda^{*}}^{-1} =
\Theta_{\lambda^{*} \hat{\lambda}}$ and $\Theta_{\lambda
\hat{\lambda}} = \Theta_{\lambda \lambda^{*}} \circ
\Theta_{\lambda^{*} \hat{\lambda}} $ this implies that
$T_{\alpha}^{\lambda} = \Theta_{\lambda
\hat{\lambda}}(T_\alpha^{\hat{\lambda}})$. Thus property (1) of
Theorem~\ref{thm:TorusWideHole1} holds.

Now we must show that if $\alpha$ is a simple closed curve in $S$ and is homotopic to the boundary, then
$T_{\alpha}^{\lambda}$ also satisfies property (1) of
Theorem~\ref{thm:TorusWideHole1}. The main step is to prove the following lemma.

\begin{lemma}
\label{lem:ExtraCaseTorusWithWideHole}
 For any two ideal triangulations $\lambda , \hat{\lambda}$ of $S$,  $
\Theta_{\lambda
\hat{\lambda}}(T_\alpha^{\hat{\lambda}})=T_\alpha^\lambda$.

\end{lemma}
\begin{proof}
The main step in the proof is the following.

\begin{lemma}
\label{lem:TorusWideHoleFirst2}
 Given two ideal triangulations $\lambda$ and $\hat{\lambda}$ which differ by only one diagonal
exchange,  $ \Theta_{\lambda
\hat{\lambda}}(T_\alpha^{\hat{\lambda}})=T_\alpha^\lambda$.
\end{lemma}

\begin{proof}
We must separate this into four cases, each of which are represented in Figure ~\ref{Fig:DiagExchangeTorusCase5}. In each case we have
$\alpha \in \Sigma(S)$ homotopic to the boundary and $\lambda$ and $\hat{\lambda}$ differing
only by exchanging edge $\lambda_i$ and $\hat{\lambda}_i$. The
edges of the quadrilaterals involved in the diagonal exchange are
$\lambda_i$, $\lambda_j$, $\lambda_k$, $\lambda_l$, $\lambda_m$,
$\hat{\lambda}_i$, $\hat{\lambda}_j$, $\hat{\lambda}_k$,
$\hat{\lambda}_l$ and $\hat{\lambda}_m$ with the quadrilaterals
edges in $\lambda$ and $\hat{\lambda}$ ordered as represented in Figure~\ref{Fig:DiagExchangeTorus2}. This
lemma follows from simple calculations nearly identical to those in Lemma ~\ref{lem:TorusWideHoleFirst} which we omit for the sake of brevity.

\end{proof}

\vspace{2mm}

\begin{figure}[htb]
\centering
 \SetLabels
 ( 0.08 * 1.02) $\lambda_j$\\
  ( .34 * 1.02) $\hat{\lambda}_j$\\
  ( 0.08 * 0.63) $\lambda_l$\\
  ( .34 * .62) $\hat{\lambda}_l$\\
( .04 * 0.795) $\lambda_i$\\
  ( .31 * .845) $\hat{\lambda}_i$\\
  ( -0.02 * 0.83) $\lambda_m$\\
  ( 0.255 * 0.83) $\hat{\lambda}_m$\\
  ( .17 * 0.83) $\lambda_k$\\
  ( .44 * .83) $\hat{\lambda}_k$\\
  ( .11 * 0.78) $\alpha$\\
  ( .39 * .86) $\alpha$\\
  (0.21 * 0.54) Case 1\\
  ( 0.66 * 1.02) $\lambda_j$\\
  ( .92 * 1.02) $\hat{\lambda}_j$\\
  ( 0.66 * 0.63) $\lambda_l$\\
  ( .92 * .62) $\hat{\lambda}_l$\\
( .615 * 0.785) $\lambda_i$\\
  ( .945 * .735) $\hat{\lambda}_i$\\
  ( 0.56 * 0.83) $\lambda_m$\\
  ( 0.835 * 0.83) $\hat{\lambda}_m$\\
  ( .755 * 0.83) $\lambda_k$\\
  ( 1.02 * .83) $\hat{\lambda}_k$\\
  ( .63 * 0.89) $\alpha$\\
  ( .95 * .89) $\alpha$\\
  (0.79 * 0.54) Case 2\\
  ( 0.08 * .33) $\lambda_j$\\
  ( .34 * 0.33) $\hat{\lambda}_j$\\
  ( 0.08 * -0.06) $\lambda_l$\\
  ( .34 * -0.06) $\hat{\lambda}_l$\\
( .13 * 0.17) $\lambda_i$\\
  ( .395 * .09) $\hat{\lambda}_i$\\
  ( -0.02 * 0.14) $\lambda_m$\\
  ( 0.255 * 0.14) $\hat{\lambda}_m$\\
  ( .17 * 0.14) $\lambda_k$\\
  ( .44 * .14) $\hat{\lambda}_k$\\
  ( .07 * 0.2) $\alpha$\\
  ( .34 * .24) $\alpha$\\
  (0.21 * -0.1) Case 3\\
  ( 0.66 * 0.33) $\lambda_j$\\
  ( 0.92 * 0.33) $\hat{\lambda}_j$\\
  ( 0.66 * -0.06) $\lambda_l$\\
  ( .92 * -.07) $\hat{\lambda}_l$\\
( .68 * 0.24) $\lambda_i$\\
  ( .905 * .235) $\hat{\lambda}_i$\\
  ( 0.56 * 0.14) $\lambda_m$\\
  ( 0.835 * 0.14) $\hat{\lambda}_m$\\
  ( .755 * 0.14) $\lambda_k$\\
  ( 1.02 * .14) $\hat{\lambda}_k$\\
  ( .62 * 0.12) $\alpha$\\
  ( .91 * .105) $\alpha$\\
  (0.79 * -0.1) Case 4\\
\endSetLabels
\centerline{\AffixLabels{
\includegraphics{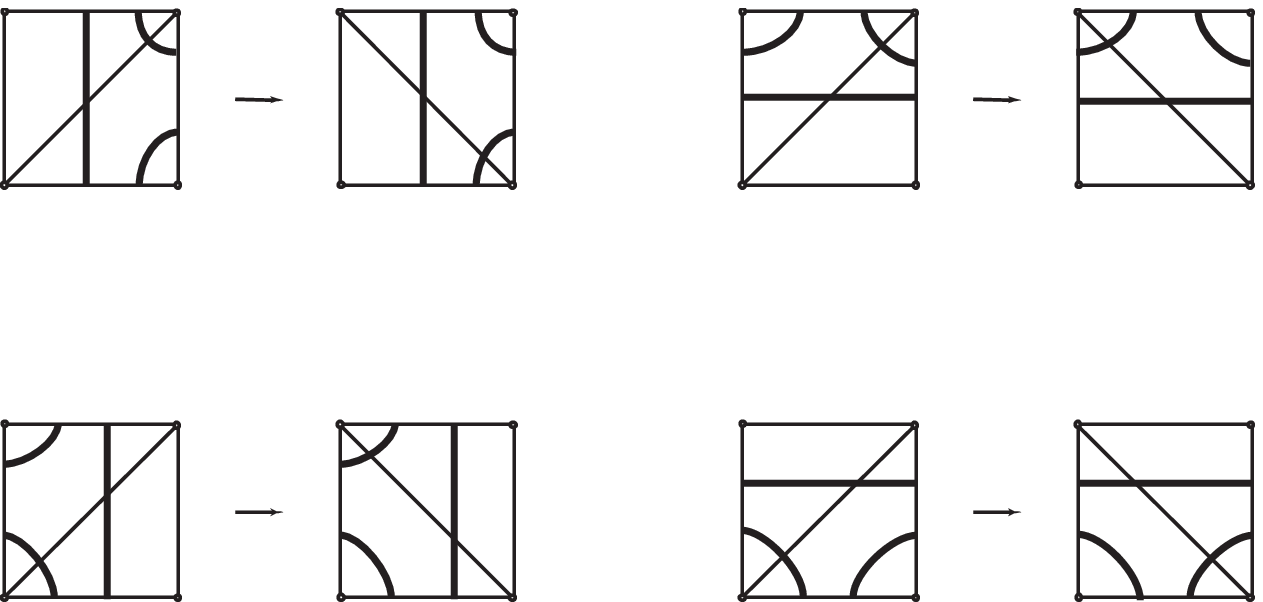}}}
\caption{} \label{Fig:DiagExchangeTorusCase5}
\end{figure}

\begin{figure}[htb]
\centering
 \SetLabels
 ( 0.2 * 1.06) $\lambda_j$\\
  ( .8 * 1.06) $\hat{\lambda}_j$\\
  ( .2 * 0.06) $\lambda_l$\\
  ( .8 * .08) $\hat{\lambda}_l$\\
( .16 * 0.56) $\lambda_i$\\
  ( .84 * .56) $\hat{\lambda}_i$\\
  ( -0.04 * 0.5) $\lambda_m$\\
  ( 1.04 * 0.5) $\hat{\lambda}_k$\\
  ( .4 * 0.5) $\lambda_k$\\
  ( .62 * .5) $\hat{\lambda}_m$\\
   ( .11 * 0.75) $T_1$\\
  ( .25 * .3) $T_2$\\
   ( .89 * 0.75) $T_1$\\
  ( .75 * .3) $T_2$\\
\endSetLabels
\centerline{\AffixLabels{
\includegraphics{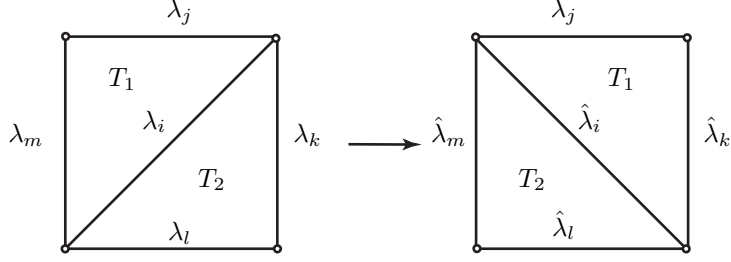}}}
\caption{Diagonal Exchange} \label{Fig:DiagExchangeTorus2}
\end{figure}

Lemma~\ref{lem:ExtraCaseTorusWithWideHole} is now a direct corollary of Theorem~\ref{thm:Penner1} and Lemma~\ref{lem:TorusWideHoleFirst2}.
\end{proof}

Properties (2) and (3) of Theorem~\ref{thm:TorusWideHole1} both
follow from definitions. Property (4) of
Theorem~\ref{thm:TorusWideHole1} follows from the way we defined
the non-quantum traces and the fact that this definition is
well-defined. This concludes the proof of
Theorem~\ref{thm:TorusWideHole1} (Theorem~\ref{thm:TorusWideHole}).
\end{proof}

We now prove Theorem~\ref{thm:TorusWideHole2} which we restate for
the sake of the reader.

\begin{figure}[htb]
\centering
 \SetLabels
 ( 0 * .2) $\alpha$\\
  ( .2 * .2) $\beta$\\
  ( .5 * .2) $\alpha\beta$\\
  ( .9 * .06) $\beta\alpha$\\
\endSetLabels
\centerline{\AffixLabels{
\includegraphics{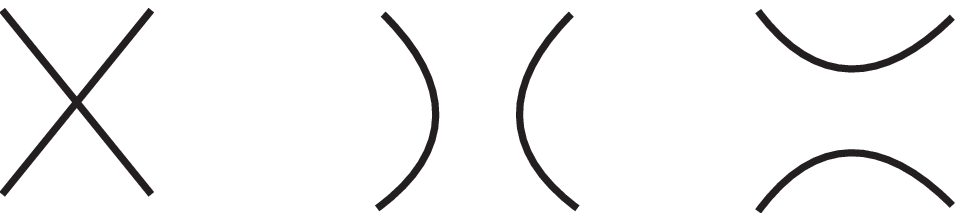}}}
\caption{Resolving Crossing in the Torus} \label{fig:SkeinTorus1}
\end{figure}

\begin{theorem}
\label{thm:TorusWideHole3} The traces $T_{\alpha}^{\lambda}$ of
Theorem~\ref{thm:TorusWideHole1} satisfy the following property: If
$\alpha$ and $\beta$ meet in one point, and if
    $\alpha \beta$ and $\beta \alpha$ are obtained by resolving
    the intersection point as in Figure~\ref{fig:SkeinTorus1}, then
    $$T_{\alpha}^{\lambda}
    T_{\beta}^{\lambda} = q^{\frac{1}{2}} T_{\alpha \beta}^{\lambda} + q^{-\frac{1}{2}} T_{\beta \alpha}^{\lambda}.$$

    In addition,
$T_{\alpha}^{\lambda}$ with $\alpha$ non-separating is the only
one which satisfy this property and conditions (1) and (4) of
Theorem~\ref{thm:TorusWideHole1}.
\end{theorem}

\begin{proof}

\begin{figure}[htb]
\centering
 \SetLabels
( 0.6 * 1.02) $\lambda_j$\\
  ( 1.05 * .5) $\lambda_k$\\
  ( 0.88 * 0.25) $\lambda_i$\\
  ( 0.72 * .15) $\lambda_h$\\
( 0.45 * .51) $\alpha$\\
  ( 0.59 * .63) $\beta$\\
  ( 0.35 * .3) $\alpha \beta$\\
  ( 0.7 * .7) $\alpha \beta$\\
  ( .3 * .67) $\beta \alpha$\\
   ( .82 * .38) $\beta \alpha$\\
\endSetLabels
\centerline{\AffixLabels{
\includegraphics{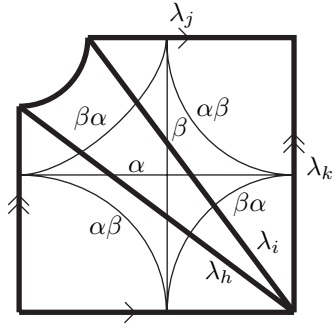}}}
\caption{Quadrilateral Labeling} \label{fig:TorusSecondThmPic}
\end{figure}

Notice that Property (1) from
Theorem~\ref{thm:TorusWideHole1} implies that it suffices to show
$T_{\alpha}^{\lambda}
    T_{\beta}^{\lambda} = q^{\frac{1}{2}} T_{\alpha \beta}^{\lambda} + q^{-\frac{1}{2}} T_{\beta \alpha}^{\lambda}$ for one
    particular $\lambda$. So, we can choose  $\lambda$ to be the ideal triangulation represented in Figure~\ref{fig:TorusSecondThmPic}. In particular, $\alpha$, $\beta$ and $\alpha \beta$ are $\lambda$-simple and $\beta \alpha$ is not.    From Theorem~\ref{thm:TorusWideHole1} $T_{\alpha}^{\lambda}$, $T_{\beta}^{\lambda}$, and $T_{\alpha \beta}^{\lambda}$ are uniquely determined and:
    \begin{equation*}
\begin{split}
T_{\alpha}^{\lambda} & =  q^{\frac{1}{4}}Z_{k,1}Z_hZ_iZ_{k,3} + q^{\frac{1}{4}}Z_{k,1}Z_h^{-1}Z_iZ_{k,3}+ q^{\frac{1}{4}}Z_{k,1}Z_h^{-1}Z_i^{-1}Z_{k,3} + q^{\frac{1}{4}}Z_{k,1}^{-1}Z_h^{-1}Z_i^{-1}Z_{k,3}^{-1} \\
T_{\beta}^{\lambda} & = q^{\frac{1}{4}}Z_{j,1}Z_hZ_iZ_{j,3} + q^{\frac{1}{4}}Z_{j,1}^{-1}Z_hZ_iZ_{j,3}^{-1} + q^{\frac{1}{4}}Z_{j,1}^{-1}Z_h^{-1}Z_iZ_{j,3}^{-1} + q^{\frac{1}{4}}Z_{j,1}^{-1}Z_h^{-1}Z_i^{-1}Z_{j,3}^{-1}\\
T_{\alpha \beta}^{\lambda} & =
Z_{k,1}Z_h^{-1}Z_i^{-1}Z_{k,3}Z_{j,1}Z_hZ_iZ_{j,3} +
Z_{k,1}^{-1}Z_h^{-1}Z_i^{-1}Z_{k,3}^{-1}Z_{j,1}Z_hZ_iZ_{j,3}
+   Z_{k,1}^{-1}Z_h^{-1}Z_i^{-1}Z_{k,3}^{-1}Z_{j,1}^{-1}Z_hZ_iZ_{j,3}^{-1}
\end{split}
\end{equation*}

\begin{figure}[htb]
\centering
 \SetLabels
( 0.13 * .82) $\lambda_i$\\
  ( .47 * .83) $\lambda_i'$\\
  ( .83 * .82) $\lambda_i''$\\
  ( 0.07 * .66) $\lambda_h$\\
  ( .59 * .22) $\lambda_h'$\\
  ( .78 * 0.2) $\lambda_h''$\\
\endSetLabels
\centerline{\AffixLabels{
\includegraphics{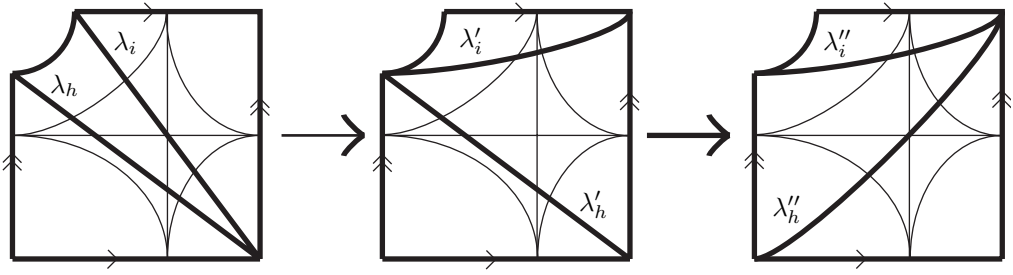}}}
\caption{Diagonal Moves in the Torus}
\label{fig:TorusSecondThmPic1}
\end{figure}

 To find $T_{\beta \alpha}^{\lambda}$ we perform two diagonal exchanges, changing $\lambda$ to
    $\lambda'$ and then $\lambda''$ respectively, as represented in Figure~\ref{fig:TorusSecondThmPic1}. Because $\beta \alpha$
    is $\lambda''$-simple $T_{\beta \alpha}^{\lambda''}$ can be determined. To calculate $T_{\beta \alpha}^{\lambda}$ we simply use
    the coordinate change maps and set $T_{\beta \alpha}^{\lambda}= \Theta_{\lambda \lambda''}(T_{\beta \alpha}^{\lambda''})$. The computations use Lemma~\ref{Lem:ComputationForCoordChangeMaps} four times, and yield:
    \begin{equation*}
\begin{split}
T_{\beta \alpha}^{\lambda} & =\Theta_{\lambda \lambda''}(T_{\beta \alpha}^{\lambda''}) =
Z_{k,1}Z_hZ_iZ_{k,3}Z_{j,1}Z_hZ_iZ_{j,3} + Z_{k,1}Z_hZ_iZ_{k,3}Z_{j,1}^{-1}Z_hZ_iZ_{j,3}^{-1} \\
&  +
Z_{k,1}Z_hZ_iZ_{k,3}Z_{j,1}^{-1}Z_h^{-1}Z_iZ_{j,3}^{-1}
+  Z_{k,1}Z_hZ_iZ_{k,3}Z_{j,1}^{-1}Z_h^{-1}Z_i^{-1}Z_{j,3}^{-1} +
Z_{k,1}Z_h^{-1}Z_iZ_{k,3}Z_{j,1}Z_hZ_iZ_{j,3}\\
&
 +   Z_{k,1}Z_h^{-1}Z_iZ_{k,3}Z_{j,1}^{-1}Z_hZ_iZ_{j,3}^{-1} +
Z_{k,1}Z_h^{-1}Z_iZ_{k,3}Z_{j,1}^{-1}Z_h^{-1}Z_iZ_{j,3}^{-1}
+  Z_{k,1}Z_h^{-1}Z_iZ_{k,3}Z_{j,1}^{-1}Z_h^{-1}Z_i^{-1}Z_{j,3}^{-1} \\
&
+  Z_{k,1}Z_h^{-1}Z_i^{-1}Z_{k,3}Z_{j,1}^{-1}Z_hZ_iZ_{j,3}^{-1} +
Z_{k,1}Z_h^{-1}Z_i^{-1}Z_{k,3}Z_{j,1}^{-1}Z_h^{-1}Z_iZ_{j,3}^{-1}
+  Z_{k,1}Z_h^{-1}Z_i^{-1}Z_{k,3}Z_{j,1}^{-1}Z_h^{-1}Z_i^{-1}Z_{j,3}^{-1}\\
&  +
Z_{k,1}^{-1}Z_h^{-1}Z_i^{-1}Z_{k,3}^{-1}Z_{j,1}^{-1}Z_h^{-1}Z_iZ_{j,3}^{-1}
+  Z_{k,1}^{-1}Z_h^{-1}Z_i^{-1}Z_{k,3}^{-1}Z_{j,1}^{-1}Z_h^{-1}Z_i^{-1}Z_{j,3}^{-1}
\end{split}
\end{equation*}

Finally we
    directly compute that $T_{\alpha}^{\lambda}
    T_{\beta}^{\lambda} = q^{\frac{1}{2}} T_{\alpha \beta}^{\lambda} + q^{-\frac{1}{2}} T_{\beta \alpha}^{\lambda}$ in the following way:

\begin{equation*}
\begin{split}
T_{\alpha}^{\lambda} T_{\beta}^{\lambda} & = q^{\frac{1}{2}} Z_{k,1}Z_hZ_iZ_{k,3}Z_{j,1}Z_hZ_iZ_{j,3} +  q^{\frac{1}{2}} Z_{k,1}Z_hZ_iZ_{k,3}Z_{j,1}^{-1}Z_hZ_iZ_{j,3}^{-1} \\
& \qquad + q^{\frac{1}{2}}
Z_{k,1}Z_hZ_iZ_{k,3}Z_{j,1}^{-1}Z_h^{-1}Z_iZ_{j,3}^{-1}
+ q^{\frac{1}{2}} Z_{k,1}Z_hZ_iZ_{k,3}Z_{j,1}^{-1}Z_h^{-1}Z_i^{-1}Z_{j,3}^{-1}\\
& \qquad + q^{\frac{1}{2}}
Z_{k,1}Z_h^{-1}Z_iZ_{k,3}Z_{j,1}Z_hZ_iZ_{j,3}
 +  q^{\frac{1}{2}} Z_{k,1}Z_h^{-1}Z_iZ_{k,3}Z_{j,1}^{-1}Z_hZ_iZ_{j,3}^{-1} \\
& \qquad + q^{\frac{1}{2}}
Z_{k,1}Z_h^{-1}Z_iZ_{k,3}Z_{j,1}^{-1}Z_h^{-1}Z_iZ_{j,3}^{-1}
+ q^{\frac{1}{2}} Z_{k,1}Z_h^{-1}Z_iZ_{k,3}Z_{j,1}^{-1}Z_h^{-1}Z_i^{-1}Z_{j,3}^{-1}\\
& \qquad + q^{\frac{1}{2}}
Z_{k,1}Z_h^{-1}Z_i^{-1}Z_{k,3}Z_{j,1}Z_hZ_iZ_{j,3}
+  q^{\frac{1}{2}} Z_{k,1}Z_h^{-1}Z_i^{-1}Z_{k,3}Z_{j,1}^{-1}Z_hZ_iZ_{j,3}^{-1} \\
& \qquad + q^{\frac{1}{2}}
Z_{k,1}Z_h^{-1}Z_i^{-1}Z_{k,3}Z_{j,1}^{-1}Z_h^{-1}Z_iZ_{j,3}^{-1}
+ q^{\frac{1}{2}} Z_{k,1}Z_h^{-1}Z_i^{-1}Z_{k,3}Z_{j,1}^{-1}Z_h^{-1}Z_i^{-1}Z_{j,3}^{-1}\\
& \qquad + q^{\frac{1}{2}}
Z_{k,1}^{-1}Z_h^{-1}Z_i^{-1}Z_{k,3}^{-1}Z_{j,1}Z_hZ_iZ_{j,3}
+  q^{\frac{1}{2}} Z_{k,1}^{-1}Z_h^{-1}Z_i^{-1}Z_{k,3}^{-1}Z_{j,1}^{-1}Z_hZ_iZ_{j,3}^{-1} \\
& \qquad + q^{\frac{1}{2}}
Z_{k,1}^{-1}Z_h^{-1}Z_i^{-1}Z_{k,3}^{-1}Z_{j,1}^{-1}Z_h^{-1}Z_iZ_{j,3}^{-1}
+ q^{\frac{1}{2}} Z_{k,1}^{-1}Z_h^{-1}Z_i^{-1}Z_{k,3}^{-1}Z_{j,1}^{-1}Z_h^{-1}Z_i^{-1}Z_{j,3}^{-1}\\
& = q^{\frac{1}{2}}T_{\alpha \beta}^{\lambda} + q^{-\frac{1}{2}}
T_{\beta \alpha}^{\lambda}
\end{split}
\end{equation*}

For any non-separating simple closed curve $\alpha$ in  $S$ the
uniqueness of $T_{\alpha}^{\lambda}$ follows from the following two facts. The
first is that property (1) from
Theorem~\ref{thm:TorusWideHole1} implies
$T_{\alpha}^{\lambda} = \Theta_{\lambda
\hat{\lambda}}(T^{\hat{\lambda}}_{\alpha})$ if $\alpha$ is
$\hat{\lambda}$-simple. The second is that Property (4) from
Theorem~\ref{thm:TorusWideHole1} implies
$T^{\hat{\lambda}}_{\alpha}$ is uniquely determined which implies
 $T_{\alpha}^{\lambda}$ is uniquely determined. This concludes
the proof of Theorem~\ref{thm:TorusWideHole3} (Theorem~\ref{thm:TorusWideHole2})
\end{proof}

\section{Spheres With Four Holes}

\begin{figure}[htb]
\centering
 \SetLabels
\endSetLabels
\centerline{\AffixLabels{
\includegraphics{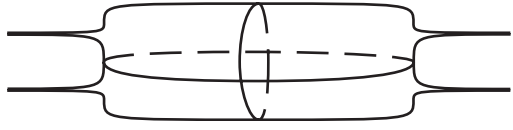}}}
\caption{The Four Times Punctured Sphere}
\end{figure}

\begin{figure}[htb]
\centering
 \SetLabels
\endSetLabels
\centerline{\AffixLabels{
\includegraphics{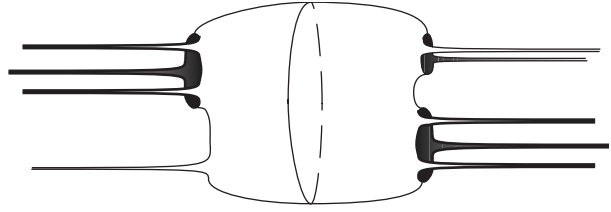}}}
\caption{A Sphere With Four Holes}
\end{figure}

Let $S$ be the surface obtained from the compact surface of genus
zero with $k$ boundary components by removing $p$ points from its
interior and at least one point from each boundary component, with
$k+p=4$, we call this surface a \emph{sphere with four holes}.
Let $\Sigma(S)$ be the set of simple closed unoriented curves in
$S$ not homotopic to the boundary.

We now state the two main theorems of the section.

\begin{theorem}
\label{thm:SphereTraces}
 Let $S$ be a sphere with four holes. There exists a family of $T_\alpha^\lambda \in \mathcal{T}_\lambda^{q^{\frac{1}{4}}}(\alpha)$, with $\lambda$
ranging over all ideal triangulations of $S$ and $\alpha$ over all
non-separating simple closed curves of $S$, which satisfies:

\begin{enumerate}
    \item If $\alpha$ is in $\Sigma(S)$ and $\lambda$ and
$\hat{\lambda}$ are two triangulations of $S$, then $
\Theta_{\lambda
\hat{\lambda}}(T_\alpha^{\hat{\lambda}})=T_\alpha^\lambda.$
    \item as $q \rightarrow 1$, $T_{\alpha}^{\lambda}$ converges to the non-quantum
    trace function $T_{\alpha}$ in $\mathcal{T}(S).$

    \item If $\alpha$ and $\beta$ are disjoint, $T_{\alpha}^{\lambda}$ and
    $T_{\beta}^{\lambda}$ commute.
    \item If $\alpha$ meets each edge of $\lambda$ at most once then $T_{\alpha}^q$
    is obtained from the classical trace $T_{\alpha}$ of Section 2.3 by multiplying each
    monomial by the Weyl ordering coefficient.

\end{enumerate}

\end{theorem}

\begin{figure}[htb]
\centering
 \SetLabels
 ( .11 * .85) $\alpha$\\
  ( .11 * .1) $\beta$\\
  ( .38 * .2) $\alpha\beta$\\
  ( .66 * .2) $\beta\alpha$\\
  ( .9 * .15) $\gamma_1$\\
  ( .9 * .47) $\gamma_2$\\
  ( .9 * .8) $\gamma_3$\\
  ( 1.02 * .5) $\gamma_4$\\
\endSetLabels
\centerline{\AffixLabels{
\includegraphics{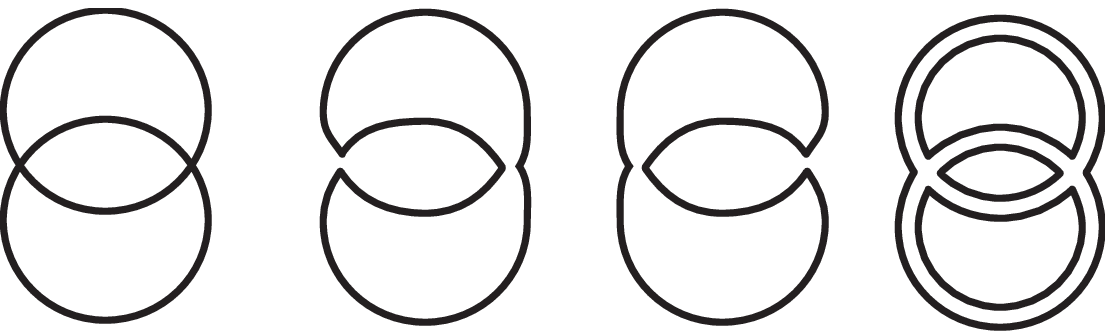}}}
\caption{Resolving Crossings in the Sphere}
\label{fig:SkeinSphere1}
\end{figure}

\begin{theorem}
\label{thm:SphereTraces2} The traces $T_{\alpha}^{\lambda}$ of
Theorem~\ref{thm:SphereTraces} satisfy the following property:

If $\alpha$ and $\beta$ meet in two points, and if
    $\alpha \beta$ and $\beta \alpha$ are obtained by resolving
    the intersection point as depicted in Figure~\ref{fig:SkeinSphere1}, then $$T_\alpha^\lambda T_\beta ^
\lambda = q T_{\alpha \beta}^\lambda + q^{-1} T_{\beta
\alpha}^\lambda +
T^{\lambda}_{\gamma_1}T^{\lambda}_{\gamma_2}+T^{\lambda}_{\gamma_3}T^{\lambda}_{\gamma_4},$$
for all $\lambda$.

In addition, $T_{\alpha}^{\lambda}$, with $\alpha$ non-separating,  is the only one which
satisfy this property and conditions (1) and (4) of
Theorem~\ref{thm:SphereTraces}.

\end{theorem}

A key result used to prove Theorem~\ref{thm:SphereTraces} is the
following proposition.

\begin{proposition}
\label{prop:KeyPropertySphere} For a curve $\alpha$ in $\Sigma(S)$
and $\alpha$-simple ideal triangulations $\lambda$ and
$\hat{\lambda}$ of $S$, there exists a sequence of $\alpha$-simple ideal
triangulations $\lambda=\lambda^0$ ,$\lambda^1$ ,$\lambda^2$ , ...
,$\lambda^{m-1}$ ,$\lambda^m =\hat{\lambda}$ such that each
$\lambda^{i+1}$ is obtained from $\lambda^i$ by a single diagonal
exchange.

\end{proposition}

\begin{figure}[htb]
\centering
 \SetLabels
 ( .5 * 1.04) $S_1$\\
  ( 1.06 * .5) $S_2$\\
  ( .5 * -0.08) $S_3$\\
  ( -0.06 * .5) $S_4$\\
  ( -0.06 * .85) $S_5$\\
\endSetLabels
\centerline{\AffixLabels{
\includegraphics{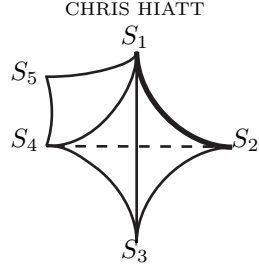}}}
\caption{Quadrilateral in the Surface} \label{fig:QuadSphere}
\end{figure}

\begin{proof}
To prove Proposition~\ref{prop:KeyPropertySphere} we must prove
the following lemma.

\begin{lemma}
\label{lem:KeyPropertySphereLemma}
    If $\lambda$ and $\hat{\lambda}$ are $\alpha$-simple ideal triangulations of $S$, there exists a sequence of $\alpha$-simple ideal triangulations
    $\lambda=\lambda^0,
    \lambda^1$ ,..., $\lambda^m$ such that $\lambda^{l+1}$ is obtained from $\lambda^l$ by a single diagonal exchange and a sequence of $\alpha$-simple ideal triangulations
    $\hat{\lambda}=\hat{\lambda}^0$ ,
    $\hat{\lambda}^1$ ,..., $\hat{\lambda}^n$ such that $\hat{\lambda}^{l+1}$ is obtained from $\hat{\lambda}^l$ by
    a single diagonal exchange and four edges $K$ and $\hat{K}$ of $\lambda^m$ and
$\hat{\lambda}^n$ respectively such that one component $C_1$ of
$S-K$ and one component $\hat{C}_1$ of
    $S-\hat{K}$ are both spheres with four spikes at infinity. Also $\lambda^m$ and $\hat{\lambda}^n$ coincide outside of $C_1$.
\end{lemma}

\begin{proof}
     We will use the following proof for both $\lambda$ and $\hat{\lambda}$.
    Consider a quadrilateral as represented in Figure~\ref{fig:QuadSphere}
    with endpoints $S_1,S_2,S_3,S_4$ where $S_1,S_2,S_3,S_4$ are all points at infinity and the edge from $S_1$ to
    $S_2$ is a boundary curve and the edges in this quadrilateral
    are edges $S_1$ to $S_2$, $S_1$ to $S_3$, $S_1$ to $S_4$, $S_2$ to $S_3$ and $S_3$ to
    $S_4$. If $S_1=S_2$ then we are done. Assume from now on that
    $S_1 \neq S_2$

    \textbf{Case 1:}  If $S_3=S_1$ and $S_4=S_1$ then doing a diagonal
    exchange on both $\lambda$ and $\hat{\lambda}$ in this quadrilateral lowers the number of edges
    ending at $S_1$ by one. Also since $\alpha$ cannot cross the edge connecting $S_1$ to $S_2$
    when you do this diagonal exchange the resulting triangulation remains $\alpha$-simple.

    \textbf{Case 2:}  If $S_1=S_3$ and $S_1 \neq S_4$, then doing a diagonal
    exchange on both $\lambda$ and $\hat{\lambda}$ in this quadrilateral lowers the number of edges
    ending at $S_1$ by two. Also, the resulting triangulation remains $\alpha$-simple.

    \textbf{Case 3:} If  $S_1 \neq S_3$ and $S_1 \neq S_4$, then doing a diagonal exchange on both $\lambda$ and $\hat{\lambda}$ in
    this quadrilateral decreases the number of edges ending at
    $S_1$ by one. Also, the resulting triangulation remains $\alpha$-simple.

    \textbf{Case 4:} If  $S_1 \neq S_3$ and $S_1=S_4$. then after doing a diagonal exchange on both $\lambda$ and $\hat{\lambda}$ in this
    quadrilateral if we consider the new quadrilateral created
    with edges $S_1$ to $S_2$, $S_2$ to $S_4$, $S_1$ to $S_4$ and
    a new point $S_5$ and edges $S_1$ to $S_5$ and $S_4$ to $S_5$
    then we see that we are again in Case 2 or Case 3. Thus after another
    diagonal exchange on both $\lambda$ and $\hat{\lambda}$ in this new quadrilateral we reduce the
    number of edges ending at $S_1$ by one or two. For the same reason as above after the first diagonal exchange the
    resulting triangulation remains $\alpha$-simple and similarly after the second diagonal exchange the resulting
    triangulation remains $\alpha$-simple.

    Now if we repeat this process until there are only two edges going to $S_1$ then we can
    effectively ``forget'' about the point at infinity $S_1$ and then repeat this process for another
    point at infinity. By the method of the proof we automatically get $\lambda^m$ and $\hat{\lambda}^n$ to coincide outside of $C_1$.

If we continue this process for each of the four wide holes of $S$
then we are left with spheres with four spikes at infinity. By the
method of the proof we automatically get $\lambda^m$ and
$\hat{\lambda}^n$ to coincide outside of $C_1$.

\end{proof}

\begin{lemma}
\label{lem:TriangleSphere} After changing the ideal triangulation
$\lambda$ by $\alpha$-simple diagonal exchanges, we can arrange so
that the following hold:
\begin{enumerate}
\item  $\lambda$ satisfies the conditions of
Lemma~\ref{lem:KeyPropertySphereLemma}, namely there exists four
edges $K$ such that one component $C_1$ of $S-K$ is a sphere with
four spikes at infinity as represented in
Figure~\ref{fig:ProofByPictureSphere}. \item  Only two edges of
$\lambda$ cross $\alpha$ \item  $\alpha$ splits $C_1$ into two
components, each of which contains four edges that are disjoint
from $\alpha$ and two boundary edges.
\end{enumerate}

\end{lemma}

\begin{figure}[htb]
\centering
 \SetLabels
( .47 * 0.87) $\alpha$ \\
\endSetLabels
\centerline{\AffixLabels{
\includegraphics{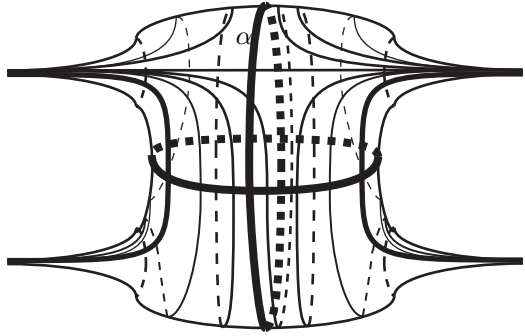}}}
\caption{Triangulation of the Sphere}
\label{fig:ProofByPictureSphere}
\end{figure}

\begin{proof}

For the first condition we simply apply
Lemma~\ref{lem:KeyPropertySphereLemma} to $\lambda$.

The second condition can be realized using moves similar to those represented in
Figure~\ref{fig:SphereReduceCrossings}. The third condition
follows from the second condition.
\end{proof}

\begin{figure}[htb]
\centering
 \SetLabels
( .08 * 0.83) $\alpha$\\
( .27 * 0.25) $\alpha$\\
\endSetLabels
\centerline{\AffixLabels{
\includegraphics{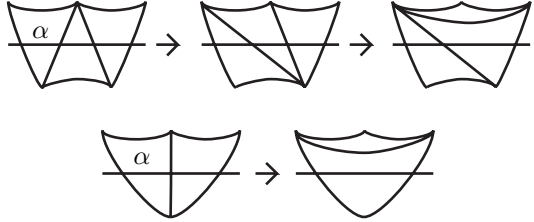}}}
\caption{Reducing Crossings for the Sphere}
\label{fig:SphereReduceCrossings}
\end{figure}

Applying Lemma~\ref{lem:KeyPropertySphereLemma}
and~\ref{lem:TriangleSphere}, we can assume without loss of
generality that $\lambda$ and $\hat{\lambda}$ satisfy the
conclusions of Lemma~\ref{lem:TriangleSphere}. By inspection, the
two edges of $\lambda$ that cross $\alpha$ must go to a single
boundary component of $C_1$ on one side of $\alpha$, and to a
single boundary component of $C_1$ on the other side of $\alpha$.

\begin{figure}[htb]
\centering
 \SetLabels
\endSetLabels
\centerline{\AffixLabels{
\includegraphics{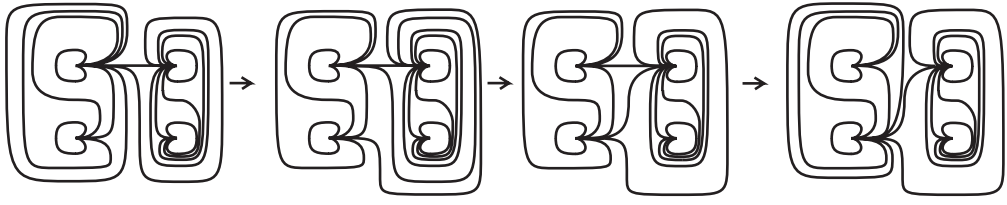}}}
\caption{Changing Spikes in the Sphere}
\label{fig:SphereChangeSides}
\end{figure}

Using the moves illustrated in Figure~\ref{fig:SphereChangeSides}
we can arrange that the edges $A$ and $B$ of $\lambda$ and
$\hat{\lambda}$ crossing $\alpha$ go to the same boundary
components of $C_1$.

\begin{figure}[htb]
\centering
 \SetLabels
\endSetLabels
\centerline{\AffixLabels{
\includegraphics{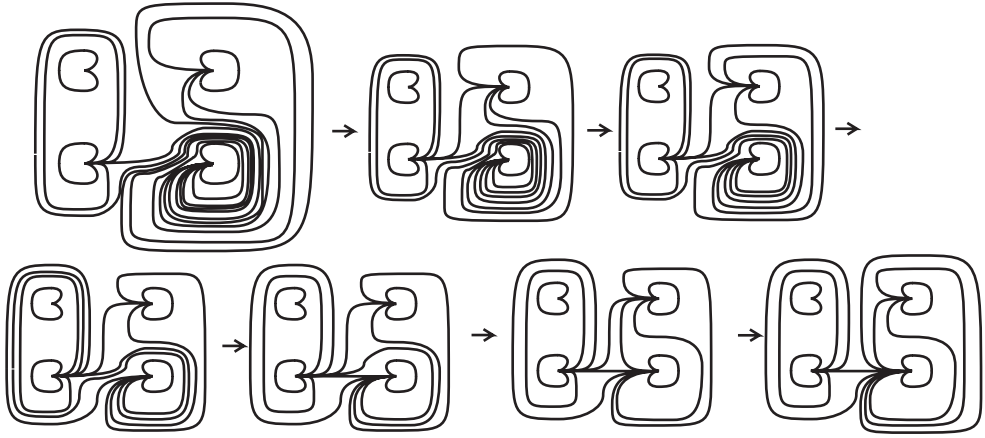}}}
\caption{Changing Windings in the Sphere}
\label{fig:ProofByPictureSphere2Colored2}
\end{figure}

\begin{figure}[htb]
\centering
 \SetLabels
\endSetLabels
\centerline{\AffixLabels{
\includegraphics{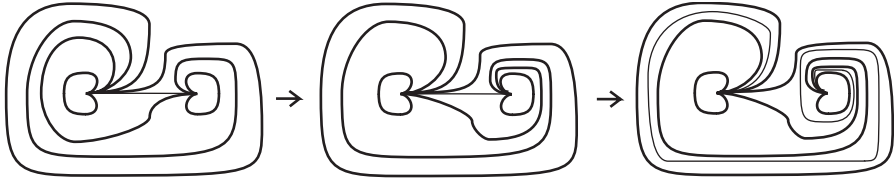}}}
\caption{Diagonal Exchange Moves in the Sphere}
\label{fig:ProofByPictureSphereKeyProperty}
\end{figure}

Finally, by using the moves represented in
Figure~\ref{fig:ProofByPictureSphere2Colored2}, we can arrange
that the two edges $A$ and $B$  wrap around $\alpha$ the same
number of times. An application of the moves illustrated in
Figure~\ref{fig:ProofByPictureSphereKeyProperty} ensures that $A$
and $B$ wrap around the boundary components of $C_1$ the same
number of times.

After these moves the two ideal triangulations now coincide.

All of this argument also works for the surface $S$ that is
obtained from the compact surface of genus zero with $k$ boundary
components by removing $p$ points from its interior and at least
one point from each boundary component, with $k+p=4$. In fact,
some of the arguments become simpler.

\end{proof}

We now prove Theorem~\ref{thm:SphereTraces} which we restate for
the sake of the reader.

\begin{theorem}
\label{thm:SphereTraces1}
 Let $S$ be a sphere with four holes. There exists a family of $T_\alpha^\lambda \in \mathcal{T}_\lambda^{q^{\frac{1}{4}}}(\alpha)$, with $\lambda$
ranging over all ideal triangulations of $S$ and $\alpha$ over all
non-separating simple closed curves of $S$, which satisfies:

\begin{enumerate}
    \item If $\alpha$ is in $\Sigma(S)$ and $\lambda$ and
$\hat{\lambda}$ are two triangulations of $S$, then $
\Theta_{\lambda
\hat{\lambda}}(T_\alpha^{\hat{\lambda}})=T_\alpha^\lambda$.
    \item as $q \rightarrow 1$, $T_{\alpha}^{\lambda}$ converges to the non-quantum
    trace function $T_{\alpha}$ in $\mathcal{T}(S)$

    \item If $\alpha$ and $\beta$ are disjoint, $T_{\alpha}^{\lambda}$ and
    $T_{\beta}^{\lambda}$ commute.
    \item If $\alpha$ meets each edge of $\lambda$ at most once then $T_{\alpha}^q$
    is obtained from the classical trace $T_{\alpha}$ of Section 2.3 by multiplying each
    monomial by the Weyl ordering coefficient.

\end{enumerate}

\end{theorem}

\begin{figure}[htb]
\centering
 \SetLabels
 ( .11 * .85) $\alpha$\\
  ( .11 * .1) $\beta$\\
  ( .38 * .2) $\alpha\beta$\\
  ( .66 * .2) $\beta\alpha$\\
  ( .9 * .15) $\gamma_1$\\
  ( .9 * .47) $\gamma_2$\\
  ( .9 * .8) $\gamma_3$\\
  ( 1.02 * .5) $\gamma_4$\\
\endSetLabels
\centerline{\AffixLabels{
\includegraphics{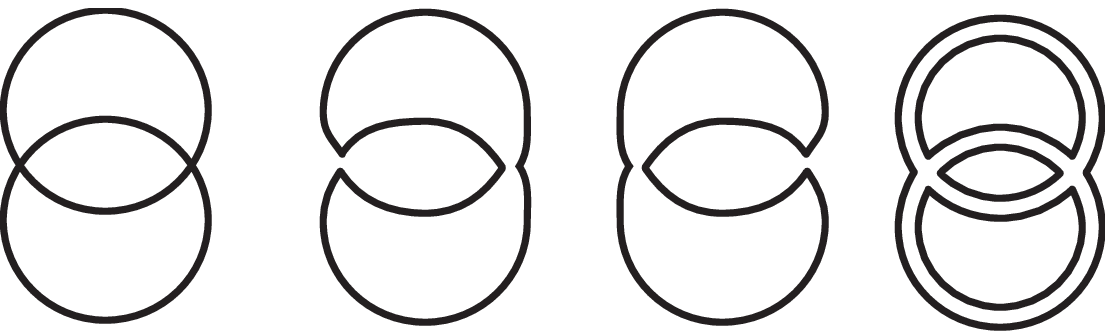}}}
\caption{Resolving Crossings in the Sphere}
\label{fig:SkeinSphere}
\end{figure}

\begin{proof}
The proof is identical to the proof of
Theorem~\ref{thm:TorusWideHole1} replacing
Proposition~\ref{prop:KeyPropertyTorus} with
Proposition~\ref{prop:KeyPropertySphere}.
\end{proof}


We conclude with the proof of Theorem~\ref{thm:SphereTraces2},
which we repeat here for the convenience of the reader.

\begin{theorem}
\label{thm:SphereTraces3} The traces $T_{\alpha}^{\lambda}$ of
Theorem~\ref{thm:SphereTraces1} satisfy the following property:

If $\alpha$ and $\beta$ meet in two points, and if
    $\alpha \beta$ and $\beta \alpha$ are obtained by resolving
    the intersection point as in Figure~\ref{fig:SkeinSphere}, then $$T_\alpha^\lambda T_\beta ^
\lambda = q T_{\alpha \beta}^\lambda + q^{-1} T_{\beta
\alpha}^\lambda +
T^{\lambda}_{\gamma_1}T^{\lambda}_{\gamma_2}+T^{\lambda}_{\gamma_3}T^{\lambda}_{\gamma_4}$$,
for all $\lambda$

In addition, $T_{\alpha}^{\lambda}$, with $\alpha$ non-separating,  is the only one which
satisfy this property and conditions (1) and (4) of
Theorem~\ref{thm:SphereTraces1}.

\end{theorem}

\begin{figure}[htb]
\centering
 \SetLabels
( .35 * 0.47) $\alpha$\\
  ( .46 * .85) $\beta$\\
\endSetLabels
\centerline{\AffixLabels{
\includegraphics{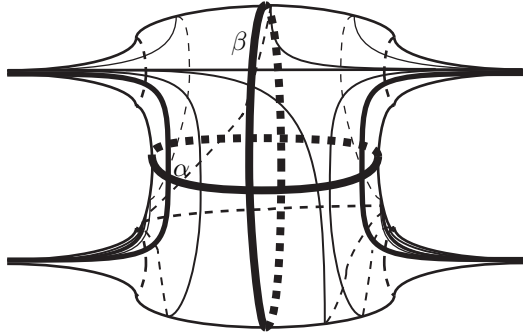}}}
\caption{$\alpha$-simple, $\beta$-simple, $\alpha \beta$-simple,
but not $\beta \alpha$-simple Triangulation}
\label{fig:ProofByPictureSphere2}
\end{figure}

\begin{figure}[htb]
\centering
 \SetLabels
( .35 * 0.47) $\alpha$\\
  ( .46 * .85) $\beta$\\
\endSetLabels
\centerline{\AffixLabels{
\includegraphics{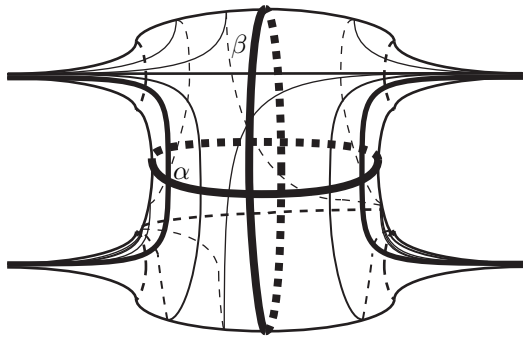}}}
\caption{$\beta \alpha$-simple Triangulation}
\label{fig:ProofByPictureSphere6}
\end{figure}

\begin{proof}

Property (1) from Theorem~\ref{thm:SphereTraces1} implies that it
suffices to show $T_\alpha^\lambda T_\beta ^ \lambda = q T_{\alpha
\beta}^\lambda + q^{-1} T_{\beta \alpha}^\lambda +
T^{\lambda}_{\gamma_1}T^{\lambda}_{\gamma_2}+T^{\lambda}_{\gamma_3}T^{\lambda}_{\gamma_4}$
for one particular $\lambda$. Let $\lambda$ be the ideal
triangulation as illustrated in
Figure~\ref{fig:ProofByPictureSphere2}. The triangulation is
$\alpha$-simple, $\beta$-simple, $\alpha \beta$-simple, but not
$\beta \alpha$-simple.

The quantum traces $T_\alpha^\lambda$, $T_\beta ^\lambda$ and
$T_{\alpha \beta}^\lambda$ are determined by Condition 4 of
Theorem~\ref{thm:SphereTraces1}. To compute $T_{\beta
\alpha}^\lambda$, we use the triangulation $\hat{\lambda}$ of
Figure~\ref{fig:ProofByPictureSphere6} and use Condition 4 of
Theorem~\ref{thm:SphereTraces1}
 to determine $T_{\beta\alpha}^{\hat{\lambda}}$ and compute $T_{\beta\alpha}^{\lambda} =
 \Theta_{\hat{\lambda}\lambda}(T_{\beta\alpha}^{\hat{\lambda}})$

At this point checking the relation $T_{\beta\alpha}^{\lambda} =
\Theta_{\hat{\lambda}\lambda}(T_{\beta\alpha}^{\hat{\lambda}})$
unfortunately requires considering $476$ terms. This computation
was verified using \textit{Mathematica}.

\end{proof}

\end{document}